\documentclass[a4paper, preprint, review, 11pt]{article}

\usepackage{booktabs}

\usepackage[linesnumbered,ruled,vlined]{algorithm2e}

\SetCommentSty{mycommfont}

\SetKwInput{KwInput}{Input}                
\SetKwInput{KwOutput}{Output}              

\usepackage{caption}
\usepackage{amsmath,amssymb}

\usepackage{graphicx}
\usepackage{adjustbox}
 \usepackage{bbm}
 \usepackage{todonotes}
 \usepackage{graphicx}
\usepackage{color}
\definecolor{cornell-red}{RGB}{179,27,27}
\usepackage{subcaption}
\usepackage{mathtools}
\usepackage{dsfont}
\usepackage{float}
\usepackage{natbib}
    \usepackage{multirow}
 
\usepackage[left=0.5in,right=0.5in,top=1in,bottom=1in]{geometry}
\usepackage{rotating}
\usepackage{xcolor}
\usepackage[colorlinks = true,linkcolor = blue,urlcolor  = blue,citecolor = blue,anchorcolor = blue]{hyperref}
\usepackage[flushleft]{threeparttable}

\usepackage{booktabs}
\usepackage{multirow}
\usepackage{rotating,tabularx}

\def\mathcolor#1#{\@mathcolor{#1}}
\def\@mathcolor#1#2#3{%
  \protect\leavevmode
  \begingroup
    \color#1{#2}#3%
  \endgroup
}
\usepackage{pifont}
\usepackage{graphicx}
\usepackage{epsfig}
\usepackage{natbib}
\usepackage{lscape}

 \bibpunct[, ]{(}{)}{,}{a}{}{,}%
 \def\bibfont{\small}%
\begin{document}




\title{Matheuristic for Vehicle Routing Problem with Multiple Synchronization Constraints and Variable Service Time}


\author{Faisal Alkaabneh$^1$ and Rabiatu Bonku$^2$}
\date{\vspace{-5ex}}

\maketitle

\author{$^1$Industrial \& Systems Engineering, North Carolina A\&T State University, Greensboro, NC 27401, United States, Email: fmalkaabneh@ncat.edu}

\author{$^1$Industrial \& Systems Engineering, North Carolina A\&T State University, Greensboro, NC 27401, United States, Email: rbonku@aggies.ncat.edu}

\abstract{%
This paper considers an extension of the vehicle routing problem with synchronization constraints and introduces the vehicle routing problem with multiple synchronization constraints and variable service time. This important problem is motivated by a real-world problem faced by one of the largest agricultural companies in the world providing precision agriculture services to their clients who are farmers and growers. The solution to this problem impacts the performance of farm spraying operations and can help design policies to improve spraying operations in large-scale farming. We propose a Mixed Integer Programming (MIP) model for this challenging problem, along with problem-specific valid inequalities. A three-phase powerful matheuristic is proposed to solve large instances enhanced with a novel local search method. We conduct extensive numerical analysis using realistic data. Results show that our matheuristic is fast and efficient in terms of solution quality and computational time compared to the state-of-the-art MIP solver. Using real-world data, we demonstrate the importance of considering an optimization approach to solve the problem, showing that the policy implemented in practice overestimates the costs by 15-20\%. Finally, we compare and contrast the impact of various decision-maker preferences on several key performance metrics by comparing different mathematical model.
}%

\noindent



%


\section{Introduction}
\label{sec:intro}
Advanced technology, especially smart farming \cite{rejeb2022drones}  and precision agriculture \cite{khanna2019evolution}, may play a significant role in making agriculture more sustainable without reducing production or farmer's revenue. According to the World Economic Forum, implementing precision agriculture may improve world production by 15–25\% by 2030 while reducing greenhouse gas emissions by 10\%. Due to this, agricultural companies are seeking innovative ways to improve their operations through precision agriculture. Agribusinesses are actively searching for methods to create new techniques or practices that will enhance the efficiency of spraying operations in large-scale farming. Operations research models have been applied to agriculture since the late 1940s and the literature is now rich with operations research application \cite{higgins2010challenges}. One of the most well-known models in operations research that finds applications across multiple fields is the Vehicle Routing Problem (VRP). An evolving form of the VRP, generally known as the vehicle routing problem with multiple synchronization constraints (VRPMS), can be used to model real-world problems involving complex requirements. The VRPMS can be applied in a wide range of applications across sectors to enhance operations efficiency, for example in logistics (see \cite{hojabri2018large}), healthcare (see \cite{hashemi2020vehicle}), service industry (see \cite{ha2020new}), ... etc, and the agriculture sector is no exception. In this work, we show how the VRPMS can be utilized to improve the efficiency of spraying operations for a company implementing precision agriculture.\\

VRPMS is a type of VRP in which the completion of a work may require the utilization of multiple vehicles \cite{drexl2012synchronization}. VRPMS models have significant practical importance due to the fact that a large variety of changes in operating rules and limitations observed in real-world applications may be modeled as VRPMS; therefore, VRPMSs have received huge scientific attention thus far. Several unresolved realistic challenges in operations research may be solved using VRPMS. Important applications of VRPMS include truck-drone delivery (\cite{masone2022multivisit}, \cite{meng2023multi}, \cite{chen2021adaptive}), healthcare and home care delivery (\cite{kim2017drone} and \cite{hashemi2020vehicle}) among others. Solving mathematical models of VRPMS is a challenging task due to the large number of decision variables and constraints involved. We observe that studies using VRPMS models have created exact, heuristic, and metaheuristic methods to address solving large-scale practical-size instances of the problem that existing solvers cannot solve within reasonable computational time. \\

Motivated by a real-world problem faced by one of the largest worldwide companies providing precision agriculture services to farmers and growers, we study the Synchronized Sprayer Tanker Routing Problem with Variable Service Time (SSTRPVST). At a high level, the SSTRPVST can be summarized as follows: two distinct types of operators (namely, a truck and a sprayer) each with a unique task are located at a central depot (namely, at the beginning of service, the sprayer and the tanker start from a known location called the starting point). Sprayers loaded with fertilizer move from this location and visit multiple locations on a farm to spray the harmful or infested plants. A sprayer stops during operation only when it has no more (i.e., not enough) fertilizer to spray the plants or when it reaches the depot at the end of the planning horizon. The tanker, on the other hand, also travels from the same central depot to the appropriate locations on the farm to meet the sprayers and refill the sprayers' tanks. The tanker does not visit every location on the farm; meaning it only travels to locations where a sprayer stops because the sprayer needs a refill. The tanker's major task is to refill the sprayer's tank wherever and whenever it is needed. It is important to note that these two vehicles are mobile; therefore, they both need to synchronize their meetings. More specifically, the tanker's route is synced with the sprayers' route, so the two vehicles meet up to perform a refilling task at a certain time and location in their sequence of operation.\\

From a modeling standpoint, vehicle routing with drones and the SSTRPVST have certain similarities in terms of task, operation, and synchronization. One significant distinction SSTRPVST possesses is the characteristic that there is no limit on the number of locations a sprayer may visit as long as the fertilizer level is not exhausted. Additionally, in the problem we solve, the refill locations where the tanker and the sprayer can meet for refilling are not known in advance. In the real-world setting, the quantity of fertilizer need to be applied at a certain location is provided in the form of a \textit{range} rather than an exact quantity to be applied. This \textit{range} specifies the lowest quantity to be applied to ensure effective spraying and the highest quantity of fertilizer to be applied to ensure safety. By incorporating this feature, we introduce the variable service time requirement. This contributes to making the problem we solve in this study difficult compared to the cases seen in the literature regarding vehicle routing with drones.\\

Yet another feature that makes SSTRPVST more challenging to solve compared to the truck-drone delivery problems is the conflict between various performance measures. More specifically, a grower aims to spray their farm as effectively as possible by applying a large quantity of fertilizer while maintaining that the applied quantity of fertilizer is within the permissible levels to ensure effective spraying. Applying a low quantity of fertilizers at an infected location within a farm does not effectively combat harmful plants and insects from growing later on. On the other hand, with more fertilizer applied, a sprayer needs more refilling which directly implies that a tanker will be traveling to perform refilling operations. Having the tanker or the sprayer traveling across a farm is considered to be a negative consequence to avoid as the wheels of these big machines affect soil compaction and it is better to minimize the traveled distance of the tanker and sprayers to maintain good soil compaction.\\

The goal of this study is to solve a new challenging problem motivated by a real-world application faced by companies implementing precision agriculture practices for spraying operations at large-scale. To this end, this study aims to develop a Mixed Integer Programming (MIP) model to formulate the SSTRPVST and solve the developed MIP efficiently to assist agribusinesses in developing plans for their spraying operations. For this difficult problem, we provide a MIP model that incorporates real-world features of the SSTRPVST without adding assumptions that yield an unrealistic model. Along with valid inequalities that are relevant to this particular problem, we focus on developing an efficient matheuristic to handle instances of large-scale problems that arise in practice because the MIP model we develop in this paper is an NP-hard, and only small-size instances of the problem can be solved within a reasonable amount of time using a solver.\\

The contributions of our paper are as follows:
\begin{itemize}
    \item We introduce a new variant of VRP with multiple synchronization constraints motivated by a real-world problem. We formally define, model, and solve the Synchronized Sprayer Tanker Routing Problem with Variable Service Time (SSTRPVST) with a matheuristic.
    \item The resulting MIP model involves a large number of variables and becomes intractable for instances of realistic size. We, therefore, propose a three-phase powerful matheuristic framework capable of handling large-scale practical size instances. The matheuristic produces fast, consistent, and high-quality solutions in a reasonable time.
    \item We generate comprehensive computational sets of experiments utilizing randomly generated sets of instances with varying settings to assess the efficacy of our approach and obtain insights. Our results demonstrate that the proposed metaheuristic can solve large cases within a fraction of the time used by the solver.  
    \item Considering real-world data, our results indicate that while there is a trade-off between quantity of fertilizer applied and traveled distance of the tanker and the sprayer, our optimization approach achieves better results on these two performance measures compared to the policy implemented in practice. Moreover, we examine the results on two different models to account for some flexibility in the requirements.
\end{itemize}

The subsequent sections of the paper are structured as follows. In section \ref{sec:litRev} we outline an overview of related literature and highlight the uniqueness of this study. Section \ref{sec:Illustrative Example} presents a detailed description of the operation we study. A mathematical model formulation for the SSTRPVST is developed in Section \ref{sec:mathModel}. Our proposed matheuristic is described in Section \ref{sec:ALNS}. The numerical analysis are reported in Section \ref{sec:Computational}. And finally, conclusions are presented in Section \ref{sec:conclude}.

\section{Literature Review}
\label{sec:litRev}

Studies in the literature on the topic of truck-drone routing and delivery  are the closest to the Synchronized Sprayer Tanker Routing Problem with Variable Service Time (SSTRPVST) as these problems share several characteristics in terms of synchronization and decisions to be made. It is worth mentioning that some studies in the literature consider one vehicle with one or multiple drones for delivery \cite{murray2015flying,agatz2018optimization,ha2018min,yurek2018decomposition,poikonen2019branch,cavani2021exact}, and \cite{roberti2021exact}, such studies are substantially different than the SSTRPVST since multiple sprayers are involved rather than one sprayer. Hence, our literature review will focus on providing a brief review of the state-of-the-art studies in the domain of truck-drone routing and delivery domain with multiple trucks.\\

\cite{tamke2021branch} creates a novel MIP model for the Vehicle Routing Problem with Drones (VRPD) with two separate time-oriented goal functions. Their paper reinforces the linear relaxation of the developed MIP by offering additional valid inequalities based on problem attributes. The authors propose a metaheuristic to solve the problem. \cite{saleu2022parallel} modifies the parallel drone scheduling traveling salesman problem by taking into consideration multiple vehicles, with deliveries divided between a vehicle and one or more drones. \cite{cokyasar2021designing} work offers a drone delivery network concept that uses automated battery-swapping machines to increase the travel range of drones. Their problem is similar to ours in that the tanker refills the sprayers' tanks in order for the sprayers to continue their praying operations.\\

\cite{schermer2019matheuristic} investigates a VRPD where a fleet of vehicles transports a certain number of drones with the aim of establishing feasible routes and drone operations that service all customers while taking the shortest amount of time. To increase solver performance, the authors structure the VRPD as a MIP with valid inequalities. They also design a metaheuristic framework to solve large instances of the problem. \cite{wang2017vehicle} and \cite{poikonen2017vehicle} study a VRPD in which a fleet of vehicles equipped with drones can be dispatched from and picked up by trucks at the depot or any of the customer locations to deliver packages to customers, with the objective of reducing the routes' maximum completion time. \cite{poikonen2020mothership} introduces the mothership and drone routing problem with the routing of a two-vehicle tandem in which the larger vehicle, which could be a ship or an airplane and the smaller vehicle, which could be a small boat or an unmanned aerial vehicle, are referred to as the drone. The authors assume fixed locations where the drones visit and return to the mothership to refuel. \cite{kitjacharoenchai2019multiple} introduces a MIP formulation that captures the delivery concept of a truck-drone combination, as well as the idea of autonomous drones flying from delivery trucks, making deliveries, and flying to any available delivery truck nearby, with the objective of reducing the arrival time of both trucks and drones at the depot after completing the deliveries.\\

Despite the similarities between the cited work so far and our work, the sprayer-routing-assigning problem we study is different from the vehicle routing problem with drones. In the papers we cited in the previous paragraphs, in the vehicle routing problem with drones, the drone is released from the truck such that each drone can visit one location, and only in a few papers two locations are visited by the drone at most, and returns to the truck or depot. On the other hand, in our problem, there is no limit on the number of nodes visited by the sprayer as long as the level of fertilizer in the sprayer's tank is enough to serve the next node. Such a restriction has a significant impact on the structure of the decision variables defined in the mathematical model. \\

Another major difference between the problem of vehicle routing with drones that has been studied in the literature and SSTRPVST is the requirement that the tanker can meet the sprayer at any node in the farm. A number of papers restrict the meeting locations to a subset of the nodes in the network, these nodes are usually called docking stations. \cite{boysen2018scheduling} creates methods for determining truck routes along robot depots and drop-off places where robots are launched. The vehicle may restock robots at these localized depots in order to deploy more until all of its users are serviced. \cite{wang2019vehicle} presents a VRPD in which a drone may travel with a vehicle, set off from its stop to serve customers, and stop at a service hub to travel with another vehicle while ensuring that the flight range and load capacity are not exceeded. \cite{liu2019optimization} studies the real-time routing of drones for meal collection and delivery in a dynamic operating context. In contrast to the aforementioned literature, where there is a predetermined subset of places known as docking stations, the refill locations in our study are not known in advance; therefore, our problem is more challenging to solve. \cite{dayarian2020same} considers the drone resupply. In their paper, the drone resupply process may occur whenever a delivery vehicle is stopped and a drone can land on the vehicle's roof. \cite{yu2022van} develops a mixed-integer program for a two-echelon van-based robot last-mile pickup and delivery system that includes time, freight, and energy, as well as an adaptive large neighborhood search method to solve larger instances and a capacity feasibility test technique for a single route. In their work, the robots may visit locations where van access is restricted and the van waits at parking nodes to drop off and/or pick up its robot, as well as to replace and/or switch its robot's batteries if necessary.\\

Other works such as \cite{ham2018integrated,chen2021adaptive}, and \cite{pina2021traveling} assume that the drones need to go back to the depot to recharge/collect more packages. For instance, \cite{nguyen2022min} develops a new optimization model known as the min-cost Parallel Drone Scheduling Vehicle Routing Problem (PDSVRP). In their work, the vehicles complete their mission in a single journey, beginning at the depot and ending at the depot, while the drones make back-and-forth trips between the depot and the consumers, one trip for each delivery. In the study by \cite{saleu2022parallel} and \cite{dell2020matheuristic} the vehicle performs a traditional delivery tour from the depot, whereas the drones are limited to back-and-forth trips in order to minimize completion time. In contrast, there are no restrictions on where the tanker and the sprayer may meet for refilling operations in the problem we analyze; hence, refill locations are more dynamic and determined through optimization. We also find other restrictions on the kinds of operations performed by the truck and the drone. For instance, \cite{poikonen2020multi} considers a single truck and multiple drones in which the truck must travel directly between successive nodes without stopping or launching additional drones. Such an assumption is also relaxed in our problem. \\

Lastly, the studies of \cite{masone2022multivisit} and \cite{gu2022hierarchical} assume that each drone is capable of visiting multiple customers before returning back to the truck. However, each vehicle is equipped with a single drone. On the other hand, the problem we study in this work assumes that one (or any) tanker can serve any sprayer without a restriction on which sprayer might be served by which tanker. Such a generalization makes our work harder than the work of \cite{masone2022multivisit} and \cite{gu2022hierarchical}.\\

To the best of our knowledge, no prior work in the literature has addressed formulating a problem similar SSTRPVST in terms of the objective function modeling and the decision variables and restrictions found in SSTRPVST.\\

\section{Operation Description}
\label{sec:Illustrative Example}

The problem we study in this work is a new problem for the OR community. Hence, in this section, we spend a great deal of time describing the problem we study to make it clear to the reader. \\

The problem we study is motivated by a problem faced by one of the largest agriculture companies worldwide specialized in providing services to farmers and growers, due to the Nondisclosure Agreement (NDA) with the company we will not be providing exact information about the name of the company, we will call that company \textit{``Company A"}. Among the services that \textit{Company A} provides is to provide spraying services to farmers to perform. Farm spraying provides several benefits to farmers, including decreased pests and weeds, healthier crops, increased crop output, and improved crop quality. \textit{Company A} receives requests from farmers to perform spraying operations in their farms around one week ahead of the visit. High-end sprayers are very expensive equipment and the cost of a sprayer tank is approximately \$800,000. This high cost of acquisition explains why it is impractical for a farmer to buy such equipment. Instead, a farmer will contract with \textit{Company A} to receive spraying services for their farm through a contract. \textit{Company A} communicates its availability to sign contracts to farmers within a state, say Iowa, in advance and farmers will be sending their requests to receive spraying services. A farmer or a grower provides \textit{Company A} with information regarding the area of their farm, and \textit{Company A} will send a team to perform scanning operations on the land using drones to get exact information on how spraying operation shall be performed, the location of harmful plants, and the topology of the farm. In return, a farmer pays \textit{Company A} an hourly rate based on how many sprayers are rented and how long the spraying operation last. \\

\begin{figure}[h]
\centering
\includegraphics[width=12cm]{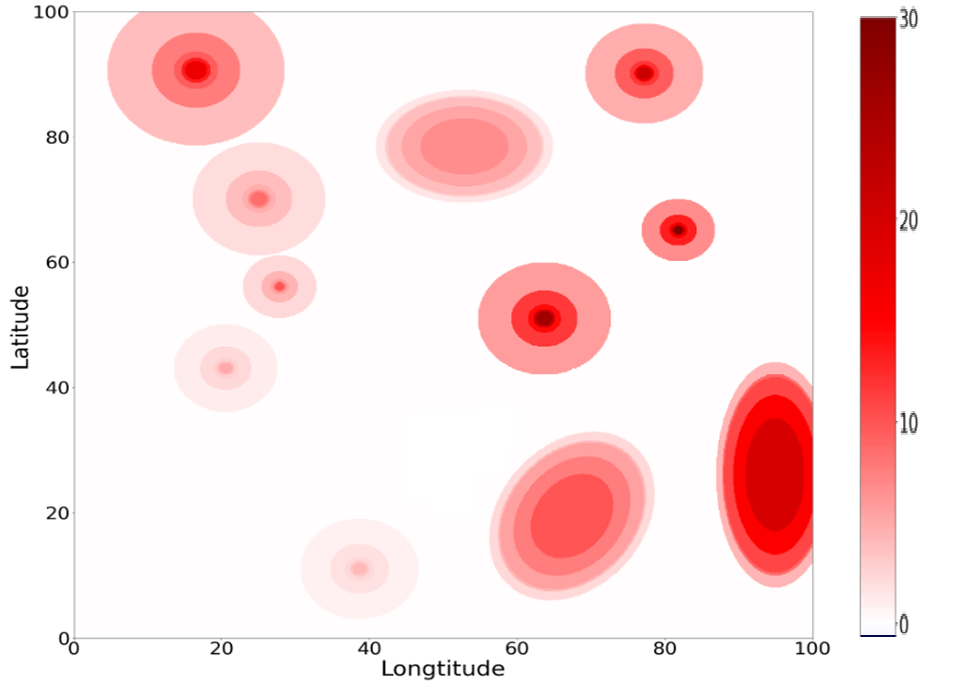}
\caption{Heat map showing infected locations at a farm}
\label{fig:heatMap}
\end{figure}

The scanning team at \textit{Company A} performs its scanning operations using drones equipped with cameras and the outcome of the scanning operations is something similar to what is depicted in Figure \ref{fig:heatMap}. As shown in Figure \ref{fig:heatMap}, there are 11 spots that are infected with harmful plants/weeds. These spots need to be sprayed to prevent the spread of diseases to the whole farm. In this work, we call each spot a \textit{location}, and later on, in the modeling section we call each location a \textit{node} for convenience. The heat map reflects the spread of weed plants and their age and therefore the estimated spraying quantity needed at each location is calculated. The spraying operations team at \textit{Company A} then estimates how much fertilizer shall be applied at that location. The spraying quantity value provided by the team is in the form of a range, namely, the team specifies how much fertilizer quantity shall be applied at least and at most. It is worth mentioning that the whole area where a weed is detected usually ends up being sprayed, this area is usually at least an acre and sometimes more. As mentioned earlier, we call each of these locations where spraying is applied a \textit{node or location}. Therefore, we represent a farm spraying operation as a graph that contains a set of \textit{nodes} that need to be served and \textit{arcs} connecting these nodes. A node is the coordinate of the center if a spot that is infested or a weed that has to be sprayed. An extensive study is needed to determine the appropriate amount of fertilizer to apply to each node using farm sprayers and the time it takes to complete the spraying operation efficiently as it impacts the overall operational cost.\\

\textit{Company A} typically uses the boom sprayer. The tank capacity of a standard boom sprayer for commercial use ranges from 1300 to 1600 US gallons. The tank is loaded with fertilizer, which is used to spray the weeds. To spray an acre of farmland, approximately 25 to 100 gallons of fertilizer are typically necessary; however, this quantity varies depending on factors such as the type of weed/crop, the depth of the infestation, and the number of diseased crops. This sprayer is mounted on a tracker with nozzles modified to provide enough coverage. When the fertilizer in the sprayer's tank runs out, a tanker comes to the node/location where a sprayer is located to refill it. Because the tanker is not a sprayer, it cannot perform any spraying operations. The tanker has a tank capacity of 4,000 to 6,000 US gallons and is loaded with the same type of fertilizer. The tanker only goes to where a sprayer needs to be refilled. The tanker's primary responsibility is to supply fertilizer to the sprayers by refilling their tanks whenever they require a chemical refill. \\

The tanker may arrive at a node/location to refill a sprayer before the sprayer performs its spraying task, during spraying, or it may arrive right on time when the sprayer's tank is empty after spraying. In these cases, the sprayer does not wait; hence the waiting time of the sprayer is zero. In case the tanker arrives after the sprayer finishes the spraying operations at a specific location, the sprayer will be waiting for a certain time. This scenario is undesirable since the sprayers are \textit{idle}, incurring additional operating costs. Following the requirement that sprayers should not be idle, in our work we impose that the waiting time of a sprayer shall be set to zero. \textcolor{red}{\textit{Company A} prefers a consistent and clear policy regarding the timing of a refill to be communicated with the sprayers, namely either a refill might take place at a node before the sprayer starts spraying or after spraying is completed at a given node. Therefore, we assume that refilling happens only after a sprayer completes spraying a node/location.} It is a challenging process to route both sprayers and tankers, as well as decide which sprayer should be allocated to which node and when and where the tanker should refill each sprayer whenever they run out of fertilizer during the operation. Overall, the SSTRPVST under consideration is quite complex, requiring decisions on routing, scheduling, and assignment to be made concurrently, as well as the interdependence between the tanker and the sprayer.\\

To provide a clear example, we use the example provided in Figure \ref{fig:heatMap} to demonstrate how to establish the assignment of sprayers to locations to be sprayed, the routing of sprayers, the refilling locations, the routing of the tanker, and finally the schedule of arrival and departure at each location for each sprayer. In this example, we assume that two sprayers are available and that there are 11 locations to be sprayed as provided in Figure \ref{fig:heatMap}.\\

\begin{figure}[!ht]
\centering
\includegraphics[width=13cm]{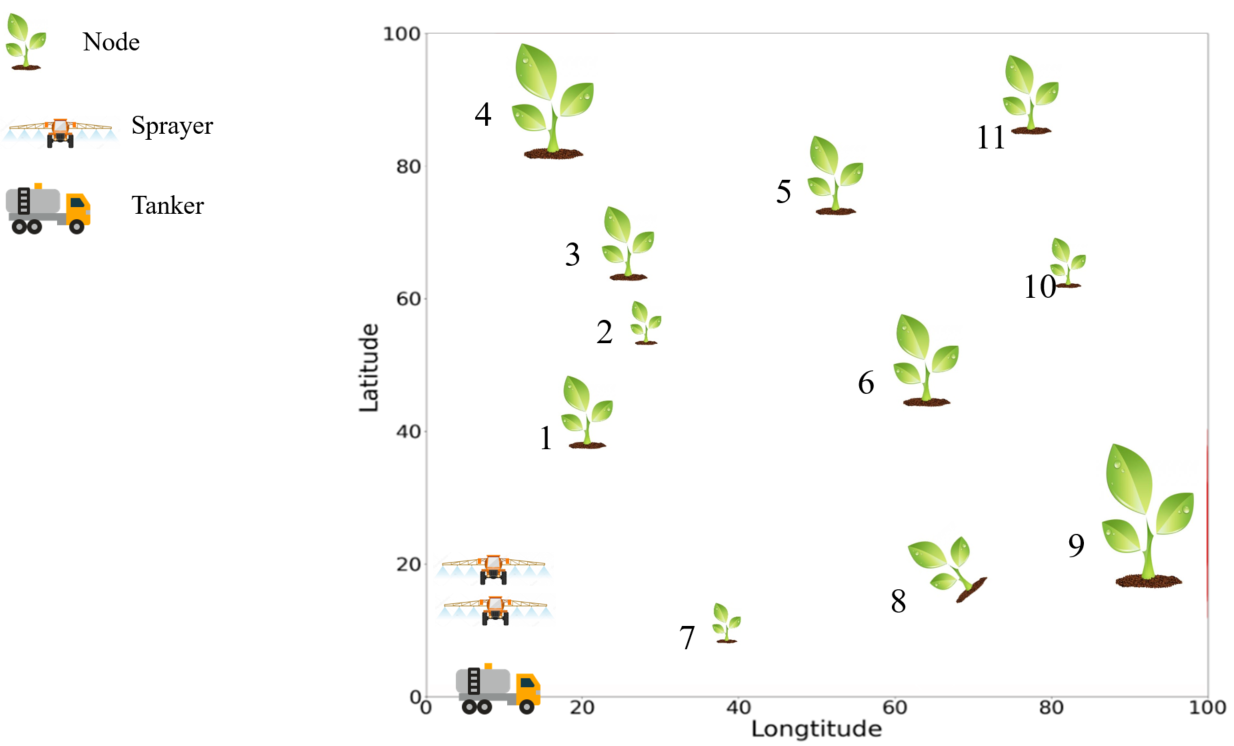}
\caption{Illustration of converting the heat map into graph}
\label{fig:step1Example}
\end{figure}

Before the start of the spraying operations, the sprayers and the tanker are located at a starting point (later on we will call that starting point \textit{the depot}). As seen in Figure \ref{fig:step1Example}, the location is simply the $(0,0)$ coordinate in the map. As shown in Figure \ref{fig:step1Example}, the area of each location to be sprayed reflects the minimum quantity of fertilizer that needs to be applied, it is hard to add more information in the figure to reflect the maximum quantity of fertilizer to be applied. Assuming that \textit{Sprayer 1} is assigned to serve locations 1, 2, 3, 4, 5, and 6 while \textit{Sprayer 2} serves locations 7, 8, 9, 10, and 11. \textit{Sprayer 1} moves to locations 1, 2, 3, and then 4,  performing spraying operations. \textit{Sprayer 1} runs out of fertilizer at location 4, thus; the tanker travels to location 4 to refill \textit{Sprayer 1's} tank. As soon as \textit{Sprayer 1's} tank is full, it resumes spraying. \textit{Sprayer 2} starts spraying at location 7 and is refilled at location 9. As of the tanker's route, the tanker leaves the starting point to location 4 to refill sprayer 1, then moves to location 9 to refill sprayer 2, and finally moves back to the starting point. Likewise, after completing the spraying operations, sprayers 1 and 2 get back to the starting point.

\begin{figure}[!ht]
\centering
\includegraphics[width=13cm]{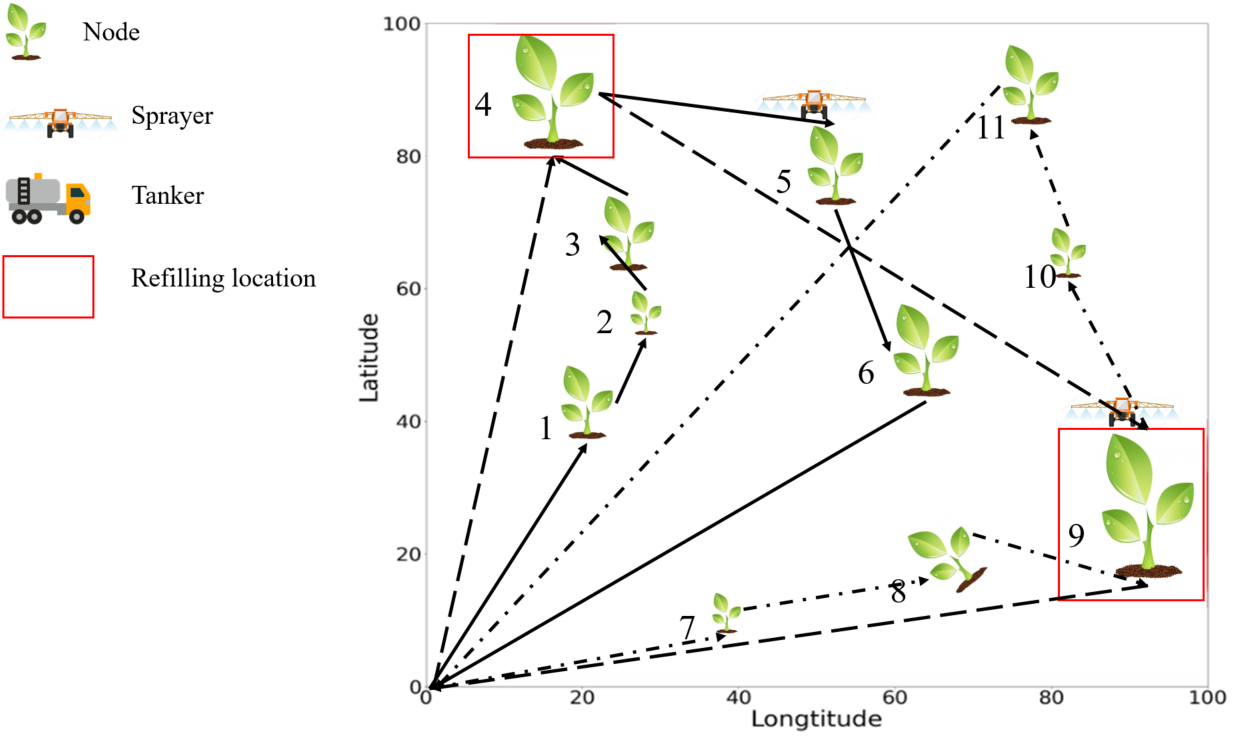}
\caption{Solution to SSTRPVST}
\label{fig:step2Example}
\end{figure}

\section{Synchronized Sprayer Tanker Routing Problem with Variable Service Time}
\label{sec:mathModel}
\subsection{Mathematical model}
\label{sec:model}
In our work, the Synchronized Sprayer Tanker Routing Problem with Variable Service Time (SSTRPVST) is defined on a complete, undirected graph $G = (\mathcal{N}, \mathcal{A})$, where $\mathcal{N}$ is the set of nodes to be sprayed in addition to the depot and $\mathcal{A}$ is the set of undirected arcs connecting these nodes. The set $\mathcal{N}=\lbrace 0, \mathcal{N}^f \rbrace$ is comprised of the depot $\lbrace 0 \rbrace$ and nodes $\mathcal{N}^f=\lbrace 1,2,\cdots,N \rbrace$. For each node to be sprayed $i\in\mathcal{N}^f$ denote the minimum quantity of fertilizer to be applied at that node as $qMin_i$ and the maximum allowable quantity of fertilizer to be applied as $qMax_i$. The calculation of $qMin_i$ and $qMax_i$ is done depending on the age and area of the weed to be sprayed as discussed in Section \ref{sec:Illustrative Example}. In this paper we assume that the spraying rate is constant, this is a mild assumption since in practice spraying rate is variable and the sprayer operator can adjust the spraying speed. We believe that ignoring the spraying rate will have a mild effect on the practical relevance of our model since the range of the spraying rate is usually not wide. We denote the constant spraying speed as $\eta^{sp}$. The working hours, which is considered as the planning horizon, is denoted by $tMax$. We assume that the network is fully connected and hence each pair of nodes in the network is connected via a direct arc. To each arc $(i, j) \in \mathcal{A}$ a travel time $t_{ij}$ is attributed. \\

There is a fleet of homogeneous sprayers and a single tanker. The set of sprayers is denoted as $\mathcal{K}$. Each sprayer is equipped with a tank and spraying equipment, the capacity of the sprayer tank is $Q^s$ gallons of fertilizer. On the other hand, a tanker's tank has a capacity of $Q^t$  gallons of fertilizer and $Q^t >> Q^s$. Finally, we define $M$ as a large number.\\

A tanker can only refill a sprayer's tank at any node after the sprayer completes serving that node. This means that a tanker can not refill a sprayer's tank while the sprayer is traveling from one node to the next. We define a binary variable $\delta_i$ that is 1 if refilling happens at node $i\in\mathcal{N}^f$ after the sprayer finishes spraying that node and 0 otherwise. \\

In reality, the time a tanker spends refilling a sprayer's tank depends on the quantity of refill. However, in our study, we assume that the time a tanker spends doing the refilling process is constant and fixed regardless of the level of fertilizer at the sprayer's tank. Since we assume a set of homogeneous sprayers, we denote the refilling time as $\xi$ and $\xi$ represents the units of time required to perform a refilling operation. In the case of a nonhomogeneous fleet of sprayers, the refilling time will have the index $k$ denoting the refilling time for each sprayer. Due to the type of engine a tanker has and a sprayer has, the tanker moves at a speed that is larger than the traveling speed of a sprayer, denote the speed factor at which the tanker is faster than a sprayer as $\beta$. In practice, a tanker has twice or more the traveling speed of a sprayer.\\

We define the following decision variables:
\begin{itemize}
    \item[-] $x_{ij}^k$ a binary routing variable equal to 1 if arc $(i, j) \in \mathcal{A}$ is traversed by sprayer $k\in \mathcal{K}$ and 0 otherwise.
    \item[-] $g_{ij}$ a binary routing variable equal to 1 if arc $(i, j) \in \mathcal{A}$ is traversed by the tanker and 0 otherwise.
    \item[-] $\delta_{i}$ a binary variable equal to 1 if refilling of  a sprayer takes place at node $i \in \mathcal{N}^f$ after serving node $i$  and 0 otherwise.
    \item[-] $\theta_{i}$ a continuous positive variable denoting the refilling start time at node $i \in \mathcal{N}^f$.
    \item[-] $y_{i}^k$ a continuous positive variable denoting the time of arrival of sprayer $k\in \mathcal{K}$ at node $i \in \mathcal{N}^f$.
    \item[-] $m_{i}^k$ a continuous positive variable denoting the waiting time of sprayer $k\in \mathcal{K}$ at node $i \in \mathcal{N}^f$.
    \item[-] $s_{i}^k$ a continuous positive variable denoting the spraying time of sprayer $k\in \mathcal{K}$ at node $i \in \mathcal{N}^f$. Note that $s_i^k$ includes the spraying activities \textit{only} and hence refilling time is not part of the spraying time.
    \item[-] $w_{i}$ a continuous positive variable denoting the time of arrival of the tanker at node $i \in \mathcal{N}^f$.
    \item[-] $h_{i}$ a continuous positive variable denoting the remaining quantity of fertilizer in the tanker's tank upon arrival to node $i \in \mathcal{N}$.
    \item[-] $l_i^k$ a continuous positive variable denoting the remaining quantity of fertilizer in sprayer's $k\in \mathcal{K}$ tank upon arrival to node $i \in \mathcal{N}^f$.
    \item[-] \textcolor{red}{$v_i$ a continuous positive variable denoting the quantity of fertilizer refilled to a sprayer whenever node $i \in \mathcal{N}^f$ is a refilling node.}
\end{itemize}

The mathematical formulation is as follows:

\begin{eqnarray}
\label{eq:Objt}
\displaystyle{\min}&& \displaystyle{\sum_{i\in \mathcal{N}}\sum_{j\in \mathcal{N}: i\neq j} t_{ij}(\sum_{k\in\mathcal{K}}x_{ij}^k+g_{ij}) + \sum_{i\in \mathcal{N}^f} \xi \delta_{i}-\sum_{k\in\mathcal{K}}\sum_{i\in\mathcal{N}^f}s_i^k}\\
s.t.
\label{eq:One}
& &\displaystyle{\sum_{k\in \mathcal{K}}\sum_{i\in \mathcal{N}:i\neq j} x_{ij}^k = 1 \quad \forall \quad j\in \mathcal{N}^f,} \\
\label{eq:Two}
& & \displaystyle{\sum_{i\in\mathcal{N}:i\neq j} x_{ij}^k = \sum_{i\in\mathcal{N}:i\neq j} x_{ji}^k \quad \forall \quad } j\in \mathcal{N},k\in \mathcal{K}, \\
\label{eq:Three}
& & \displaystyle{\sum_{i\in\mathcal{N}:i\neq j} g_{ij} = \sum_{i\in\mathcal{N}:i\neq j} g_{ji} \quad \forall \quad  j\in \mathcal{N},}\\
\label{eq:Two2}
& & \displaystyle{\sum_{j\in\mathcal{N}^f} x_{0j}^k = 1 \quad \forall \quad } k\in \mathcal{K}, \\
\label{eq:Three3}
& & \displaystyle{\sum_{j\in\mathcal{N}^f} x_{j0}^k = 1 \quad \forall \quad } k\in \mathcal{K}, \\
\label{eq:Four}
& & \displaystyle{\delta_{i} = \sum_{j\in \mathcal{N}:j\neq i} g_{ji} \quad \forall \quad i\in\mathcal{N}^f,}  \\
\label{eq:Five}
& &\displaystyle{\sum_{i\in \mathcal{N}^f} g_{0i} = \sum_{i\in \mathcal{N}^f} g_{i0}= 1,}\\
\label{eq:Six}
&&\displaystyle{y_i^k + s_i^k  \leq \theta_{i} \quad \forall i\in \mathcal{N}^f,k\in \mathcal{K},}\\
\label{eq:Seven}
& & \displaystyle{y_j^k\geq y_i^k+s_i^k+ t_{ij}x_{ij}^k + \xi \delta_{i} + m_i^k \color{red}-\color{black} M(1-x_{ij}^k) \quad \forall \quad i,j\in \mathcal{N}^f:  i\neq j,k\in \mathcal{K},}\\
\label{eq:Eleven}
& & \displaystyle{y_j^k \leq y_i^k + s_i^k + t_{ij}x_{ij}^k + \xi\color{red}\delta_i\color{black} + m_i^k + M(1-x_{ij}^k),\quad \forall \quad i,j\in\mathcal{N}^f:i\neq j,k\in\mathcal{K}}  \\
\label{eq:Twelve}
& & \displaystyle{y_i^k \leq t_{0i}x_{0i}^k + M(1-x_{0i}^k) \quad \forall \quad i\in \mathcal{N}^f,k\in\mathcal{K},}\\
\label{eq:Eight}
& & \displaystyle{y_i^k\geq t_{0i}x_{0i}^k \quad \forall \quad i\in \mathcal{N}^f,k\in \mathcal{K},}\\
\label{eq:Nine}
& & \displaystyle{\theta_{i} \geq t_{0i}/\beta g_{0i} \quad \forall \quad i\in \mathcal{N}^f,}  \\
\label{eq:Ten}
& & \displaystyle{s_i^k + y_i^k + \color{red}t_{i0}x_{i0}^k\color{black}\leq tMax  \quad \forall \quad i\in \mathcal{N}^f,k\in \mathcal{K},}\\
\label{eq:Fourteen}
& & \displaystyle{\theta_i + (t_{ij}/\beta+\xi)g_{ij} - M(1-g_{ij}) \leq w_j \quad \forall  i\in \mathcal{N}^f, \quad j\in\mathcal{N}:i\neq j}\\
\label{eq:Thirteen}
& & \displaystyle{w_i \leq \theta_i  \quad \forall \quad i\in \mathcal{N}^f,}\\
\label{eq:TwentyOne}
& & \displaystyle{\theta_{i} - (y_i^k +s_i^k) - M(1-\delta_{i})\leq m_i^k = 0\quad \forall \quad i\in \mathcal{N}^f,k\in\mathcal{K}}\\
\label{eq:Sixteen}
& & \displaystyle{\eta^{sp}s_i^k \leq qMax_i\sum_{j\in \mathcal{N}:j\neq i}x_{ji}^k\quad \forall \quad i \in\mathcal{N}^f,k\in\mathcal{K},}  \\
\label{eq:Seventeen}
& & \displaystyle{\eta^{sp}s_i^k \geq qMin_i\sum_{j\in \mathcal{N}:j\neq i}x_{ji}^k\quad \forall \quad i \in\mathcal{N}^f,k\in\mathcal{K},}  \\
\label{eq:EighteenP1}
& & \displaystyle{\color{red} l_j^k \leq l_i^k - \eta^{sp}s_i^k + v_i+Q^s(1-x_{ij}^k)\quad \forall \quad i,j \in\mathcal{N}^f:i\neq j,k\in\mathcal{K},}  \\
\label{eq:EighteenP2}
& & \displaystyle{\color{red} l_j^k \geq l_i^k - \eta^{sp}s_i^k + v_i-Q^s(1-x_{ij}^k)\quad \forall \quad i,j \in\mathcal{N}^f:i\neq j,k\in\mathcal{K},}  \\
\label{eq:refillQuant1}
& & \displaystyle{\color{red} v_i \leq -l_i^k + \eta^{sp}s_i^k + Q^s+Q^s(1-\delta_{i})\quad \forall \quad i \in\mathcal{N}^f,k\in\mathcal{K},}  \\
\label{eq:refillQuant2}
& & \displaystyle{\color{red} v_i \leq Q^s\delta_i\quad \forall \quad i \in\mathcal{N}^f,k\in\mathcal{K},}  \\
\label{eq:EighteenTanker}
& & \displaystyle{h_j \leq h_i - v_{i} + M(1-g_{ij})\quad \forall \quad i,j \in\mathcal{N}^f:i\neq j,}  \\
\label{eq:Nineteen}
& & \displaystyle{l_0^k = Q^s,}  \\
\label{eq:Twenty}
& & \displaystyle{h_0 = Q^t\quad \forall k\in\mathcal{K},}  \\
\label{eq:FTen19}
& & \displaystyle{l_j^k \leq Q^s  \quad \forall \quad j\in \mathcal{N}^f,}  \\
\label{eq:EighteenTanker2}
& & \displaystyle{h_j \geq Q^t - M(1-g_{0j})\quad \forall \quad j \in\mathcal{N}^f,}  \\
\label{eq:FIntg}
& &\displaystyle{g_{ij}, x_{ij}^k,\delta_{i} \in \lbrace 0,1\rbrace}, \theta_{i},s_{i}^k,y_{i}^k,l_{i}^k,w_{i}, h_{i}\geq 0 \quad \forall i,j,k. \
\end{eqnarray}

The objective function (\ref{eq:Objt}) minimizes the traveled distance of the sprayers and the tanker by aiming at minimizing the traveling time metric and maximizes the productivity of the sprayers by maximizing the service time of sprayers at each node. The refilling time is also minimized since refilling a sprayer requires the sprayer to stop working and hence it is considered a waste of time. \textcolor{red}{The company we partner with divides the time a sprayer spends into three parts to demonstrate the efficiency of a sprayer. These parts are: idle time (the time a sprayer spends waiting to get a refill), up-time (the time a sprayer spends spraying), and traveling time. As such, no priority or weight is given to any part over the others, hence, our model does not assign a weight of importance to any of these parts over the other}. The SSTRPVST constraints can be divided into four parts:\\

\textbf{Sprayers and tanker routing constraints (\ref{eq:One})-(\ref{eq:Five}):} Constraints (\ref{eq:One}) ensure that each node is visited by one and only one sprayer. Constraints (\ref{eq:Two}) are the inflow and outflow constraints to guarantee that a sprayer leaves a node that it visits. Likewise, constraints (\ref{eq:Three}) are the inflow and outflow constraints for the tanker. Constraints (\ref{eq:Two2}) and (\ref{eq:Three3}) guarantee all sprayers leave and return to the depot. Constraints (\ref{eq:Four}) impose the requirement that a node can be visited by the tanker if and only if that node is a refilling node. Constraints (\ref{eq:Five}) imply that there is only one tanker.\\

\textbf{Sprayers and tanker scheduling constraints (\ref{eq:Six})-(\ref{eq:TwentyOne}):} Constraints (\ref{eq:Six}) set the refilling start time to be the arrival time of sprayer $k$ to node $i$ plus the service time in the earliest. Note that this constraint is essential to set the refilling time to be performed after service completion. Constraints (\ref{eq:Six}) can be also an equality constraint since the waiting time of the sprayer shall be 0 but we leave it as inequality to accommodate for more modeling variation in case the decision maker allows waiting time for the sprayer. Constraints (\ref{eq:Seven}) and (\ref{eq:Eleven}) set the arrival time (and hence the service start time) of node $j$ that is visited by sprayer $k$ after serving node $i$. Constraints (\ref{eq:Twelve})-(\ref{eq:Eight}) set the starting service time of the first node visited by a sprayer right after departing from the depot. These constraints also imply that the traveling time from the depot to the first node is part of the planning horizon. Likewise, constraints (\ref{eq:Nine}) set the arrival time of the first node visited by the tanker. Constraints (\ref{eq:Ten}) ensure that each sprayer completes their tasks and gets back to the depot within the maximum allowable time. Constraints (\ref{eq:Fourteen}) calculate the arrival time of the timer at the next node after performing a refilling operation at the previous node. Constraints (\ref{eq:Thirteen}) ensure that refilling starting time at node $i$ is only permissible after the tanker arrives at node $i$. Constraints (\ref{eq:TwentyOne}) guarantee that the waiting time of sprayers is set to 0. \\

\textbf{Sprayers' and tanker's tank capacity and refilling constraints (\ref{eq:Sixteen})-(\ref{eq:FIntg}):} Constraints (\ref{eq:Sixteen}) and (\ref{eq:Seventeen}) ensure that quantity of fertilizer applied at node $i$ does not exceed the maximum allowable quantity or is less than the minimum quantity to be applied. Note that we do not define a new decision variable to denote the sprayed quantity at any node,  instead, we use the service time and the spraying rate to express the quantity of fertilizer applied at a node. Constraints (\ref{eq:EighteenP1}) and (\ref{eq:EighteenP2}) track the level of fertilizer of sprayer $k$ upon arriving at node $j$ after visiting node $i$, constraints (\ref{eq:EighteenP1}) and (\ref{eq:EighteenP2}) also account for the event of refilling at node $i$ prior to arriving at node $j$. Constraints (\ref{eq:refillQuant1}) calculate the quantity of refill for a sprayer at node $i$ based on the remaining level of fertilizer after serving node $i$. Constraints (\ref{eq:refillQuant2}) set the appropriate value of the refilling quantity to be at most $Q^s$ in case there is a refilling at node $i$ and 0 otherwise. Constraints (\ref{eq:EighteenTanker}) track the level of fertilizer in the tanker's tank. Constraint (\ref{eq:Nineteen}) and (\ref{eq:Twenty}) ensure that the tank of the sprayers and the tanker is full after departing from the depot. Constraints (\ref{eq:FTen19}) ensure that the sprayer's tank does not hold more fertilizer than its capacity. Constraints (\ref{eq:EighteenTanker2}) imply that a tanker leaves the depot with a full tank. The variable domains are given in constraints (\ref{eq:FIntg}).

\subsection{Valid Inequalities}
Adding valid inequalities to mathematical model (\ref{eq:Objt})-(\ref{eq:FIntg}) helps to further reduce the search space size. In this study, we add a set of valid inequalities based on an upper bound of the maximum service time by any sprayer (or alternatively, the maximum quantity of fertilizer a sprayer can perform spraying operations). While constraints (\ref{eq:Eleven}) limit the total time a sprayer travels, gets refilling, and performs spraying, stronger inequalities can be added to limit the service time of a sprayer. Let $s'$ be the maximum service time a sprayer can spend. Valid inequalities to mathematical model (\ref{eq:Objt})-(\ref{eq:FIntg}) can be written as:

\begin{equation}
    \sum_{i\in\mathcal{N}^f}\eta^{sp}s_i^k \leq  s' \quad\forall k \in \mathcal{K}
\end{equation}

The value of $s'$ is obtained by solving a problem similar to the Orienteering Problem (OP) with minor modifications to account for the refilling time. The modified OP model can be formally defined on an undirected graph $G=(\mathcal{N}_{OP},\mathcal{A})$ where $\mathcal{N}_{OP} = 0 \cup \mathcal{N}^f\cup O_{ter}$ is the set of the depot where the trip starts ($0$), nodes to potentially be visited ($\mathcal{N}^f$), and the terminal node ($O_{ter}$). Similar to the notation we use in mathematical model (\ref{eq:Objt})-(\ref{eq:FIntg}), we denote the set of nodes including the depot and all nodes to sprayed as $\mathcal{N}=0\cup\mathcal{N}^f$, clearly set $\mathcal{N}$ does not include the terminal node $O_{ter}$. Note that the terminal node is, in fact, a copy of the depot to return to after the trip. We set the reward associated with node $i\in\mathcal{N}^f$ as $qMax_i/\eta^{sp}$ whereas $qMax_0 = qMax_{O_{ter}} = 0$. We use $tMax$ as the maximum total travel time allowed to complete a tour and we use the traveling time parameter $t_{ij}$ as defined above. As for the decision variables, let $s_i$ be a positive decision variable denoting the service time at node $i\in\mathcal{N}^f$, $x_{ij}$ be a binary variable set to 1 if node $j$ is visited after visiting node $i$, $numF$ as a positive integer variable to count the number of refilling times needed, and finally define integer variable $u_i$ to track the order of visiting node $i$ to prevent subtours. The modified OP mathematical model we use is: 

\begin{eqnarray}
\label{eq:OPbj}
\displaystyle{\max}&& \displaystyle{\sum_{i\in \mathcal{N}^f}s_i}\\
s.t.
\label{eq:OP1}
& &\displaystyle{s_i \leq qMax_i/\eta^{sp}\sum_{j\in\mathcal{N}:j\neq i} x_{ji} \quad \forall i \in \mathcal{N}^f,} \\
\label{eq:OP2}
& & \displaystyle{\sum_{i\in \mathcal{N}^f} x_{0i} = \sum_{i \in \mathcal{N}^f} x_{i0_{ter}}}, \\
\label{eq:OP3}
& & \displaystyle{\sum_{i\in \mathcal{N}^f} x_{0i} = 1}, \\
\label{eq:OP4}
& & \displaystyle{\sum_{i\in 0\cup\mathcal{N}^f} x_{ij} = \sum_{i\in \mathcal{N}^f\cup 0_{ter}} x_{ji} \quad\forall j \in \mathcal{N}^f}, \\
\label{eq:OP5}
& & \displaystyle{\sum_{i \in \mathcal{N}}\sum_{j \in \color{red}\mathcal{N}_{OP}\color{black}} t_{ij}x_{ij} + \sum_{i\in \mathcal{N}^f}s_i + \xi*numF \leq tMax}, \\
\label{eq:OP6}
& & \displaystyle{\sum_{i\in \mathcal{N}^f}(\eta^{sp}s_i)/Q^s\leq numF}, \\
\label{eq:OP7}
& & \displaystyle{1\leq u_i\leq \vert \mathcal{N}^f\vert +1 \quad\forall i\in\mathcal{N}^f\cup 0_{ter}}, \\
\label{eq:OP8}
& & \displaystyle{u_i-u_j+(\vert \mathcal{N}^f\vert+1)x_{ij} \leq \vert \mathcal{N}^f\vert \quad \forall i \in \mathcal{N}^f\cup 0_{ter},\quad \forall j \in \mathcal{N}^f\cup 0_{ter}: i\neq j}, \\
\label{eq:OP10}
& &\displaystyle{x_{ij} \in \lbrace 0,1\rbrace, \quad q_{i} \geq 0 \quad u_i,numF \in \mathbf{Z}^+.} \
\end{eqnarray}

The value of the objective function of model (\ref{eq:OPbj})-(\ref{eq:OP10}) is set to be $s'$ which is an upper bound on the service time performed by any sprayer. In case the capacity of the sprayers is heterogeneous, the value of $Q^s$ should be adjusted accordingly. Constraints (\ref{eq:OP1}) ensure that the maximum quantity of applied fertilizer at node $i$ does not exceed the maximum allowable quantity of fertilizer. Constraints (\ref{eq:OP2}) imply that the tour starts by leaving the depot and ends by getting back to the dummy copy of the depot (i.e., the sprayer's tour starts and ends at the depot). Constraints (\ref{eq:OP3}) imply that the sprayer leaves the depot to visit one of the nodes in the set $\mathcal{N}^f$. Constraints (\ref{eq:OP4}) determine that each location is visited at most once and the inflow and outflow are preserved. Constraints (\ref{eq:OP5}) ensure feasibility of the maximum time a sprayer may spend before getting back to the depot. In particular, constraints (\ref{eq:OP5}) imply that sprayer's traveling time plus service time plus refilling time are within the maximum time a sprayer can spend. Constraints (\ref{eq:OP6}) calculate the number of refilling times of the sprayer's tank. Constraints (\ref{eq:OP7}) and (\ref{eq:OP8}) are the subtour elimination constraints. And finally, constraints (\ref{eq:OP10}) are the variables' domain constraints.

\subsection{Lower bound}
\label{sec:LowerBound}
Due to the complexity of mathematical model (\ref{eq:Objt})-(\ref{eq:FIntg}) and the fact that our developed matheuristic does not provide exact information on the converging to the optimal point during the search process. We opt to establish lower bounds for model (\ref{eq:Objt})-(\ref{eq:FIntg}) to measure the efficiency of our developed matheuristic. A lower bound based on linear relaxation of some of the binary variables of model (\ref{eq:Objt})-(\ref{eq:FIntg}) is expected to be poor since model (\ref{eq:Objt})-(\ref{eq:FIntg}) involves big-M. Instead, we use a lower bound based on relaxing some of the constraints of model (\ref{eq:Objt})-(\ref{eq:FIntg}). Specifically, we relax the requirement that a tanker needs to be available for refilling operations and instead we assume that once a sprayer's tank is out of fertilizer, it gets refilled instantaneously. Therefore, we remove decision variables $g_{ij},\theta_i,w_i,$ and $h_i$ and we drop the corresponding constraints needed.

\section{Matheuristic for Sprayer-Tanker Routing Problem}
\label{sec:ALNS}
Clearly, the problem we present in this study is an NP-hard problem since it is a VRP with multiple synchronization constraints. Using a MIP solver such as Gurobi enables us to solve small-scale instances of the mathematical model presented in Section \ref{sec:mathModel}. On the other hand, solving large-scale practical instances of the mathematical model calls for developing a decomposition approach or a powerful matheuristic. To this end, in this study we focus on developing an efficient matheuristic to find high-quality solutions of mathematical model (\ref{eq:Objt})-(\ref{eq:FIntg}) within reasonable computational time. \textit{Company A} runs spraying operations once for a farmer per season; and hence the SSTRPVST is a tactical problem and the farmers decide on the spraying plan a week or two in advance. Given the fact that Company A contracts with several farmers in each state, an acceptable computational time should not exceed 10 hours per farmer.  \\

The matheuristic we develop consists of three phases, the first phase aims at finding an initial feasible solution that respects all constraints (initialization), the second phase aims at improving the initial feasible solution iteratively using ALNS combined with local search, and the \textcolor{red}{improvement} stage solves a mathematical model based on the information collected during the ALNS search process (improvement).\\

\subsection{Calculating service time, refilling points, and tanker routes for fixed sprayers' routing}
\label{sec:serviceTime}
In this subsection, we introduce the method we implement to calculate the service time at each node (i.e., $s_i^k$), the refilling points (i.e., $\delta_i$), and to establish the tanker's route (i.e., $g_{ij}$). For a given set of routes of sprayers, we calculate the values of $s_i^k,\delta_i$ and $g_{ij}$ using the procedure outlined in this subsection. But first, we need to introduce some notation. 

A valid way to represent a solution of SSTRPVST is to state the route of each sprayer as an ordered sequence
$$r_k=[\bigl< 0,0,0,Q^s,0 \bigr>,\bigl< i,y_i^k,s_i^k,l_i^k,\delta_i \bigr>, \bigl< j,y_j^k,s_j^k,l_j^k,\delta_j \bigr>, ..., \bigl< 0,y_0^k,u_0,l_0^k,\delta_0 \bigr>],$$

where each tuple $\bigl< i,y_i^k,u_i,l_i^k,\delta_i \bigr>$ represents the node number $i\in\mathcal{N}^f$, starting time of service of node $i$, the quantity of fertilizer applied at node $i$, the level of the sprayer's $k$ tank upon arriving at node $i$, and if refilling is performed at node $i$, respectively. Note that the nodes in each $r_k$ are served according to their order given in the route sequence (i.e., $y_i^k < y_j^k$). 

For tanker's route representation, we use the following notation:
$$r_{tanker}=[\bigl< 0,0,0 \bigr>,\bigl< i,a_i^k,h_i \bigr>, \bigl< j,a_j^k,h_j \bigr>, ..., \bigl< 0,a_0,h_0 \bigr>],$$

where each tuple $\bigl< i,a_i^k,h_i \bigr>$ represents the node number $i\in\mathcal{N}^f$, starting time of refilling sprayer's $k$ tank at node $i$, and the level of tank of fertilizer upon arriving at node $i$, respectively. Note that the nodes in each $r_{tanker}$ are served according to their order given in the route sequence (i.e., $a_i^k < a_j^k$).

\begin{itemize}
    \item Assume that the capacity of sprayer $k$ tank is set to $\alpha Q^s$ for $0< \alpha  \leq 1$. This assumption implies that there is buffer capacity at sprayer $k$ tank equals to $(1-\alpha)Q^s$ denoted as $B_k$. For the time being, assume that the sprayer's tank capacity is $\alpha Q^s$ and its buffer capacity is $B_k=(1-\alpha)Q^s$. For any value of buffer capacity $B_k$, we can easily calculate the buffer time as $bt_k=B_k/\eta^{sp}$ \textcolor{red}{(lines \ref{l1} and \ref{l2} in Algorithm \ref{STalgo}).} 
    \item Assume that the demand of each node is set to $qMin_i$. For sprayer's $k$ route, follow the order of nodes to be visited and assume that the service time at each node is set to $qMin_i/\eta^{sp}$, a refilling is only performed at end of service at node $i$ if and only if sprayer $k$ leaves node $i$ to node $j$ such that $h_j < qMin_j$ and $x_{ij}^k=1$. For the time being, we assume that the refilling process happens instantaneously and we ignore the synchronization constraints \textcolor{red}{(lines \ref{r11}-\ref{r22} in Algorithm \ref{STalgo}).}
    \item Following the calculation in the previous step, we can calculate the \textit{\textcolor{red}{earliest} arrival time} at each node. We call that time \textit{\textcolor{red}{earliest} } since it will not be the actual arrival time after taking into account the requirement that the tanker needs to be available for the refill process to start. 
    \item Using the \textcolor{red}{earliest} arrival time calculated in the previous step, we order the refilling points based on the refilling starting time in descending order and the tanker route is established following that order \textcolor{red}{(lines \ref{lineT1}-\ref{lineT2} in Algorithm \ref{STalgo}).} 
    \item The actual refilling time can now be calculated by considering the availability of the tanker \textcolor{red}{(lines \ref{tArTime}-\ref{refillArTime}) in Algorithm \ref{STalgo})}. If the waiting time of a sprayer at any refilling point is greater than zero, we use some of the buffer time to increase the service time for the nodes that are visited prior to the refilling point and up the previous refilling point without violating constraints (\ref{eq:Seventeen}) \textcolor{red}{(lines \ref{feas1}-\ref{feas2} in Algorithm \ref{STalgo})}. On the other hand, if the waiting time is larger than the buffer time, the solution is considered infeasible \textcolor{red}{(lines \ref{infeas1}-\ref{infeas2} in Algorithm \ref{STalgo}) and the objective function calculated will include a penalty term due to violation of constraints (\ref{eq:TwentyOne}), we provide more details on handling infeasibility in Section \ref{sec:infeasibleSol}}. 
    \item \textcolor{red}{For a new value of $\alpha$, lines \ref{alphaSearch}-\ref{bestAlpha} decide if the service times and tanker route under the given value of $\alpha$ are better than the previously searched values and update the values of service times, tanker route, and objective function value.}  
\end{itemize}

Note that in the previous process, the larger the value of $\alpha$ the fewer refills will be needed and hence the less the tanker travels. On the other hand the larger the value of $\alpha$ the more service time can be performed. In order to find a good value of $\alpha$, in this work we perform a line search method to find the best value of $\alpha$. We summarize the process of calculating the service time at each node, refilling points, and tanker's route in Algorithm \ref{STalgo}. \textcolor{red}{The core function of Algorithm \ref{STalgo} is to find good values of the service time, refilling points, and tanker's route given the sprayers' routes through a line search by varying the value of $\alpha$. Later in Section \ref{sec:infeasibleSol} we present the function we use to calculate the objective function value of any solution during the execution of the ALNS.}

\begin{algorithm}[H]
  \renewcommand{\arraystretch}{0.5}
  \small
\DontPrintSemicolon
  \KwInput{Sprayers' routes}
  \KwOutput{Service time at each node, refilling points, and tanker's route}
  $objFuncVal \gets +\infty,\delta_i,s_i^k,g_{ij} \gets \emptyset$\;
  \For{$\alpha \in \lbrace 0.60,0.65,0.70,...,1.0\rbrace$}
   {       $\Bar{Q}_s \gets \alpha Q^s$\;
            \label{l1}            
            \For{$k \in \mathcal{K}$}
            { $B_k \gets (1-\alpha) Q^s , bt_k\gets (1-\alpha)Q^s$\;\tcp*{Calculate the \textit{expected} arrival time, service time, and tank level for each sprayer}
            \label{l2}
              \For{$i \in \mathcal{N}$}
            { \label{r11}
                \For{$j \in \mathcal{N}^f$}
                    {
                    \If{$x_{ij}^k == 1$}
                        { 
                            $l_j^k \gets l_i^k-qMin_i$\;
                            \If{$l_i^k-qMin_i < 0$}
                            {   
                            $l_j^k \gets \Bar{Q}_s, \Bar{\delta}_i \gets 1, y_j^k\gets \Bar{s}_i^k+t_{ij}+y_i^k+\xi$\; \tcp*{Refilling shall be performed at the previous node}
                            }
                            \Else {
                            
                            $l_j^k \gets l_i^k-qMin_i, \Bar{\delta}_i \gets 0, y_j^k\gets \Bar{s}_i^k+t_{ij}+y_i^k$\; \tcp*{No refilling needed at the previous node}
                            \label{r22}
                            }
                        }
                
            } 
            }  
            }
            \For{$i \in \mathcal{N}^f \& \Bar{\delta}_i == 1$}
                    {\label{lineT1}
                    Insert $i$ in $tankerRoute$ based on calculated $y_i^k$ \;
                
            }
            Calculate $\Bar{g}_{ij}$ based on $tankerRoute$.\; 
            \label{lineT2}
            \For{$i \in \mathcal{N}$}
                    {
                    \For{$j \in \mathcal{N}^f \& \Bar{g}_{ij} == 1$}
                    {
                    $w_j \gets \theta_i+t_{ij}+\xi$\;
                    \label{tArTime}
                    $\theta_j \gets \max\lbrace w_j,y_j^k+\Bar{s}_j^k\rbrace$\;
                    \label{refillArTime}
                    \If{$w_j > y_j^k+\Bar{s}_j^k \& w_j - y_j^k-\Bar{s}_j^k < bt_k$}
                    {
                    \label{feas1}
                    Increase the value of all $\Bar{s}_i^k$ prior to $j$ up till $\Bar{\delta}_i=1$ within the route of sprayer $k$.\;
                    Update the value of all $y_i^k$ prior to $j$ up till $\Bar{\delta}_i=1$ within the route of sprayer $k$.\;
                    $\theta_j \gets y_j^k+\Bar{s}_j^k$\;
                    $\Bar{m}_i^k \gets y_j^k+\Bar{s}_j^k-w_j$\;
                    \label{feas2}
                
                    } 
                \Else
                {\label{infeas1}
                Problem Infeasible\;
                $Penalty = \lambda \sum_{k\in\mathcal{K}}\sum_{i\in\mathcal{N}^f}\Bar{m}_i^k$\;
                    } 
                \label{infeas2}
                    }
                Calculate the objective function value\; \label{alphaSearch}
                $TempObjFuncVal \gets \sum_{i\in \mathcal{N}}\sum_{j\in \mathcal{N} i\neq j} t_{ij}(\sum_{k\in\mathcal{K}}x_{ij}^k+ \Bar{g}_{ij}) + \sum_{i\in \mathcal{N} i} \xi \delta_{i}-\sum_{k\in\mathcal{K}}\sum_{i\in\mathcal{N}^f}\Bar{s}_i^k+Penalty$\;
                \If{$TempObjFuncVal < objFuncVal$}
                    {
                    $objFuncVal \gets TempObjFuncVal $\;
                    $\delta_i\gets \Bar{\delta}_i,s_i^k\gets \Bar{s}_i^k,g_{ij} \gets \Bar{g}_{ij}$\;
                    \label{bestAlpha}
                
                    }

   }
   \;
   \label{line2}
   }
\caption{Calculating service time, refilling points, and tanker's route}\label{STalgo}
\end{algorithm}

\subsection{Initialization phase}
\label{sec:feasible solution}
Finding a good initial solution is an important step to establishing a high-quality solution of the developed matheuristic. In our work, we propose two techniques to generate an initial feasible solution. The first technique encompasses clustering the nodes into clusters based on the number of available sprayers (i.e., the size of set $\mathcal{K}$), assigning a cluster for each sprayer, performing routing operations for the sprayers, and finally establishing refilling decision variables and tanker's routing decision. The second technique is based on greedily inserting one node at a time in the best location in a sprayer route. We provide the details of each technique in this subsection. We call the first technique \textit{clustering initialization} and the second technique \textit{greedy initialization}.\\

The clustering initialization technique mimics the cluster-first-route-second technique in which initially a set of clusters are established based on $k$-clustering algorithm to group the nodes based on their location (i.e., $(x,y)$ coordinates), a sprayer is assigned one cluster, routing decisions are created based on solving a TSP instance for each cluster to generate the shortest route connecting all the nodes within a specific cluster. This procedure generates routing decisions for sprayers. In order to generate a feasible solution, the missing part is then to calculate the service time for each node, refilling positions, and tanker's routing. These calculations are made following the procedure outlined in Section \ref{sec:serviceTime}. Note that an infeasible solution might be produced out of this process, specifically, the generated solution may violate constraints (\ref{eq:Eleven}). To generate a feasible solution, we remove nodes from the busiest cluster to a cluster with less number of nodes. This is due to the fact that the $k$-clustering algorithm may produce non-balanced clusters. The second technique we implement to generate an initial feasible solution is the \textit{greedy initialization}. The greedy initialization tries to find the best position of a node.\\

At the beginning of the construction heuristic, each vehicle route $r_k$ contains only tuples for starting and ending the tour in the depot, i.e. [$\bigl< 0,0,0,Q^s,0 \bigr>$, $\bigl< 0,0,0,Q^s,0 \bigr>$]. The goal of the heuristic is to sequentially insert one node in one position of one of the $\mathcal{K}$ routes (one route for each sprayer) such that all constraints (\ref{eq:One})-(\ref{eq:FIntg}) are satisfied. Each time a new sequence is generated, we run the procedure in Section \ref{sec:serviceTime} to calculate service time, refilling locations, and the tanker's route.


\subsection{Adaptive Large Neighborhood Search}
The second phase of our developed matheuristic is based on using Adaptive Large Neighborhood Search (ALNS) metaheuristic to further improve the initial solution generated following the procedure in Section \ref{sec:feasible solution} as well as generating high-quality solutions through a search process. The ALNS metaheuristic proved to be an effective metaheuristic to solve variations of VRP, see \cite{sacramento2019adaptive,kitjacharoenchai2020two,frey2022vehicle,alkaabneh2022multi} and \cite{alkaabneh2023routing}. The ALNS was originally developed by \cite{ropke2006adaptive} and it works well for various settings of complex VRP due to the flexibility ALNS provides to design specific features for any problem.\\ 

In the first phase, an initial solution is passed to the ALNS which implements a series of iterations to find a better solution by sequentially destroying a given solution and then repairing it, this iterative process is called the improvement phase. The selection of a destroy operator to use at any iteration as well as the selection of a repair operator follows a roulette-wheel mechanism in which the probability of selecting an operator depends on its performance so far in the search process. If the new repaired solution is accepted according to an acceptance criterion, the current solution is replaced by the new solution and the procedure starts again.\\

Our ALNS incorporates innovative features in the improvement phase. Specifically, we deviate from the standard ALNS presented by \cite{ropke2006adaptive} in three ways: (1) we allow temporarily infeasible solutions, however penalty term is added in the objective function to account for the infeasibilities; (2) we develop new destroy operators, which exploit the underlying problem structure; and (3) we design an intensification search step to find more promising solutions by solving a simple mathematical model. As reported by \cite{cordeau2001unified}, the main benefit of exploring infeasible solutions is to enable better traversing of the search space because it reduces the chance of getting stuck in a local optimum and by oscillating between feasible and infeasible regions with the appropriate penalty terms. For the SSTRPVST, we allow infeasibilities due to violations of no waiting time of sprayers constraints, see constraints (\ref{eq:TwentyOne}).\\

As an overview of the organization of this section, subsections \ref{sec:Destroy} and \ref{sec:Repair} present the destroy and repair operators in our ALNS, respectively. Subsection \ref{sec:localSearch} describes the technique implemented in our work to perform local search. Subsection \ref{sec:acceptCriteria} details the acceptance criteria of our ALNS metaheuristic. Subsection \ref{sec:infeasibleSol} discusses how the developed ALNS metaheuristic handles infeasible solutions during the search process. Subsection \ref{sec:weightUp} describes the updating of the weights, and finally, we present the pseudo code in section \ref{sec:pseudo code}.

\subsubsection{Destroy operators}
\label{sec:Destroy}

A destroy operator removes a certain number of nodes from the provided solution and add them to the removal list. One advantage of the ALNS metaheuristic is the flexibility in selecting which set of nodes shall be removed from a solution and added to the removal list. The number of nodes to be removed as well as the characteristics of nodes depends on which destroy operator is selected. Each destroy operator is designed in a unique way to select a set of nodes to be removed. In this work, our developed ALNS utilizes 11 destroy operators, some of them are adapted from the literature while others are designed specifically for SSTRPVST.

\begin{figure}[h]
\centering
\includegraphics[width=16cm]{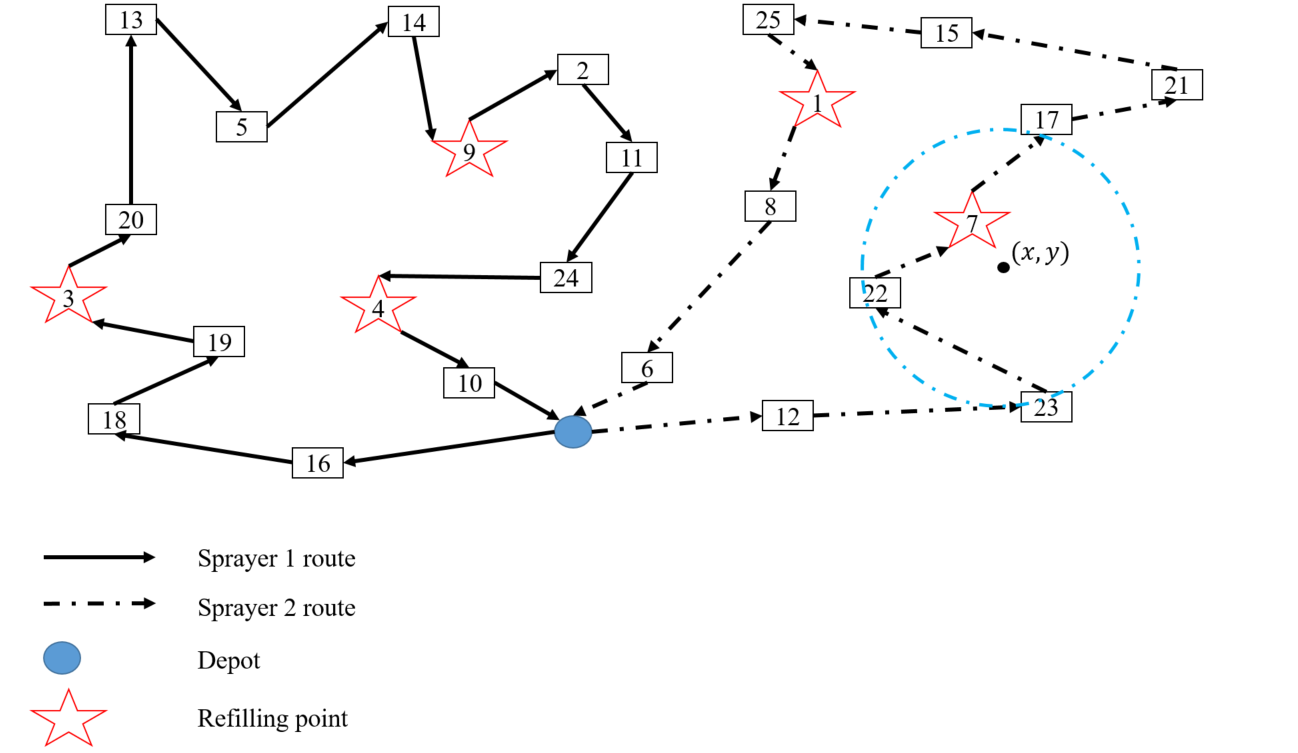}
\caption{Example to illustrate destroy operators}
\label{fig:DesOp}
\end{figure}

The operators we use from the literature are: random, route, longest distance, worst distance, historical knowledge, and zone, see \cite{ropke2006adaptive,demir2012adaptive}, and \cite{zhang2018electric}. For SSTRPVST, the main characteristic creating difficulties is the synchronization between the tanker and sprayers for refilling; therefore, in our study we introduce four operators specifically addressing this spatial and temporal complexity: refilling points, nodes proceeding a refilling point, nodes succeeding a refilling point, and $\kappa-$nodes around a refilling point. To clearly illustrate the concept of each of the new destroy operators, we use the example provided in Figure \ref{fig:DesOp}. Figure \ref{fig:DesOp} depicts an example of a feasible solution for two sprayers and 25 nodes example and how each destroy operator selects nodes from that solution.\\

\textbf{1. Random}\\
The random removal destroy operator is one of the simplest operators. Random removal operator picks nodes at random, removes them from a solution, and adds these nodes to the removal list $\mathcal{L}$. The number of nodes to be removed using a random removal operator is $U*\vert \mathcal{N}^f\vert$, where $U$ is a random number from a uniform distribution with a lower bound of $0.07$ and an upper bound of $0.15$. To illustrate the mechanism of random destroy operator using the example given in Figure \ref{fig:DesOp}, the operator will start by selecting a member of the set $\mathcal{N}^f$, suppose that node 5 was selected, remove it from the solution, add it to the removal list of $\mathcal{L}$, and the process is repeated until $p$ nodes are removed. For example, $\mathcal{L} = \lbrace 5, 11, 17, 22, 23\rbrace$.\\

\textbf{2. Route}\\
The route destroy operator picks one sprayer at random, removes all nodes served by the selected sprayer, and adds the removed nodes to the removal set $\mathcal{L}$. Using the example in Figure \ref{fig:DesOp}, the operator may select route 1 to be removed, and hence the following set of nodes will be removed from the solution and inserted in the removal list; therefore, $\mathcal{L} = \lbrace 2,3,4,5,10,11,13,14,16,18,19,20,24\rbrace$. \\

\textbf{3. Longest distance-service time}\\
The longest distance-service time destroy operator picks nodes that cause long traveling time to and from while little contribution to the objective function. More specifically, the operator calculates the distance traveled to each node and from that node to the next node minus the service time (i.e., $d_i=t_{i^-,i}+t_{i,i^+}-s_i^k$ where nodes $i^-$ and $i^+$ are the preceding and successive nodes of node $i$, respectively, according to node's $i$ current position). Then the operator sorts nodes based on $d_i$ in ascending order and removes $p$ nodes with the longest distance and short service times. This saves valuable extra space on the route, which might subsequently take up a node with a shorter distance and a higher service time. The size of the removal list $\mathcal{L}$ is controlled by the value of $p$ and $p=U*\vert \mathcal{N}^f\vert$, where $U$ is a random number from a uniform distribution with a lower bound of $0.05$ and an upper bound of $0.12$.\\

\textbf{4. Worst distance}\\
\textcolor{red}{This destroy operator removes nodes with high-cost savings after removal. The saving is defined as the total traveling time to be saved by a sprayer due to removing node $j$ due to traveling between the node $j$ and the preceding node (denote this node $j^-$) and the following node (denote this node $j^+$) on the route of a node. The operator selects and removes node $j = \max_{j\in \mathcal{N}^f}  {|t_{j^-,j} + t_{j,j^+}|}$.}\\ 

\textbf{5. Historical knowledge}\\
This operator removes nodes placed in positions that lead to bad solutions in previous iterations. At every iteration during the execution of the ALNS, the operator calculates the position cost of each node $i$ which is equal to the sum of the distances between its prior and succeeding nodes and is computed as $\pi_{i} = t_{j,i} + t_{i,e}-s_i^k$, where nodes $j$ and $e$ are the preceding and successive nodes of node $i$, respectively, according to node's $i$ current position. The operator updates the node $i$ position cost $\pi_{i}$ to be the lowest value of all $\pi_{i}$ values determined up to any iteration. The operator selects a node $j^{*}$ on a route with the maximum difference from its node cost i.e. $j^{*} = \max_{j\in \mathcal{N}^f} {\pi_{j} - \pi_{j^*}}$ and adds to the removal list $\mathcal{L}$, the procedure continues until $p$ nodes have been added to the list.\\

\textbf{6. Zone}\\
The zone removal operator adds a set of nodes to the removal list based on their proximity to a random point generated by the operator. More specifically, a random point across the map is generated $(x,y)$, and all nodes enclosed by a circle with radius $rd$ are added to the removal list. Using the example in Figure \ref{fig:DesOp}, the point $(x,y)$ is generated and nodes $7,17,22,$ and $23$ are encapsulated within the circle; therefore, these nodes are removed from the solution and added to the removal list $\mathcal{L} = \lbrace 7,17,22,23\rbrace$.\\

\textbf{7. Refill Position Removal}\\
This method is specially designed for our \textcolor{red}{SSTRPVST}. It adds all refilling points to the removal list $\mathcal{L}$. Using the example in Figure \ref{fig:DesOp}, the refilling points are nodes $1,3,4,7,$ and $9$; these nodes will be removed from the solution and added to the removal list $\mathcal{L}=\lbrace 1,3,4,7,9\rbrace$.\\

\textbf{8. Neighbors of refill positions}\\
In this operator, each refilling point, it's successor node, and it's preceding node are added to the removal list $\mathcal{L}$. Using the example in Figure \ref{fig:DesOp}, the refilling points are nodes $1,3,4,7,$ and $9$. Each of these refilling points is succeeded by a node and proceeded by a node within the sprayer routes. For example, node's $3$ predecessor is node 19 and its successor is node 20. Likewise, node 1, its predecessor is node 25 and its successor is node 8. Hence, the removal list $\mathcal{L}=\lbrace 19,3,20,14,9,2,24,4,10,22,7,17,25,1,8\rbrace$.\\

\textbf{9. $\kappa-$Neighbors of refill positions}\\
\textcolor{red}{This operator is a generalization of operator 8 by increasing the number of nodes to be removed that are within the vehicle sequence of refilling points.} This operator removes a refilling point and the $\kappa$ nodes prior to and following a refilling point. We set the value of $\kappa = 2$. Using the example in Figure \ref{fig:DesOp}, the refilling points are nodes $1,3,4,7,$ and $9$. For node 1, the two nodes proceeding it are nodes 15 and 25, and the two nodes succeeding it are nodes 8 and 6. The removal list of nodes will be $\mathcal{L}=\lbrace 23,22,7,17,21,15,25,1,8,6,18,19,3,20,13,5,14,9,2,11,24,4,10\rbrace$.\\

\textbf{10. Subtour prior to refilling position}\\
Inspired by \cite{frey2022vehicle}, we design the subtour prior to a refill position operator to pick a refilling point and all nodes visited prior to the selected node all the way up to the depot or the previous refilling point. This operator picks half of the refilling points at random and performs the insertion process of the selected nodes to the removal list $\mathcal{L}$. Using the example in Figure \ref{fig:DesOp}, the refilling points are nodes $1,3,4,7,$ and $9$. Assuming that nodes 1 and 9 are picked at random. For node 1, the nodes preceding it all the way to a refilling point or a depot are nodes 25, 15, 21, and 17. The removal list of nodes will be $\mathcal{L}=\lbrace 17,21,15,25,1,20,13,5,14,9\rbrace$. \\

\textbf{11. Subtour following a refilling position}\\
Similar to subtour prior to the refilling position operator, the subtour following to refilling position operator removes a refilling point and all nodes following that refilling node in a route up to the depot or the next refilling point. Using the example in Figure \ref{fig:DesOp}, the refilling points are nodes $1,3,4,7,$ and $9$. Assuming that nodes 4 and 7 are picked at random. For node 4, the nodes succeeding all the way to a refilling point or a depot are nodes 10. The removal list of nodes will be $\mathcal{L}=\lbrace 4,10,7,17,21,15,25\rbrace$.\\

\subsubsection{Repair Methods}
\label{sec:Repair}
In this section, we present the two repair operators utilized in our ALNS metaheuristic. Insertion/Repair operators are used to repair the partially destroyed solution by re-inserting the nodes in the removal list $\mathcal{L}$ sequentially. The two insertion operators we use were inspired by \cite{ropke2006adaptive}. See the brief description below:\\

\textbf{1. Greedy}\\
The greedy repair operator receives the removal list $\mathcal{L}$ and a destroyed solution. It iteratively selects a node from the list to insert one at a time into existing routes in the best position. The node and insertion position with the lowest insertion cost is chosen, and the corresponding insertion operation is performed.\\

\textbf{2. Regret}\\
The regret repair operator evaluates the best and second best positions to insert a node from the removal list $\mathcal{L}$ into the current solution. The node with the highest difference between the best and second best position is inserted first and the evaluation is repeated again until all nodes from the removal list $\mathcal{L}$ have been inserted.\\

The selection of which repair operator to use at any iteration is done in a random way and we do not assign weights of importance to the repair operator (i.e., the probability of selecting any operator at any iteration is always $0.5$ regardless of its previous performance). We opt to fix the weight of each repair operator since there are only two operators.

\subsubsection{Intensive Local Search}
\label{sec:localSearch}
In our ALNS to further intensify the search process we use a mathematical model to improve on some solutions that are promising. More specifically, whenever a new \textcolor{red}{best} solution is better than the current best solution found so far during the search process, we explore the neighbor of that solution to further improve the new best solution by solving mathematical model (\ref{eq:Objt})-(\ref{eq:FIntg}) with certain binary variables being fixed.\\

Finding optimal routes of the sprayers alone does not guarantee finding the optimal solution to model (\ref{eq:Objt})-(\ref{eq:FIntg}) since finding the optimal refilling nodes, the service time at each node, and the route of the tanker is an optimization problem by itself. Since establishing the refilling points and calculating the service time at each node is done using a simple rule during the search process (see discussion in Section \ref{sec:serviceTime}) it is very critical to implement a way to find optimal refilling points and calculate the service time at each node for a given sprayers' routes. In this work, we design a novel local search technique to find near-optimal refilling points and calculate the service time at each node given a set of sprayers' routes. The local search we design to find optimal service time and refilling positions for given sprayers' routing decisions is part of the intensification phase of our matheuristic.\\

The idea of our novel local search is based on the fact that eliminating nodes from the search space as candidates for refilling points speeds up solving mathematical model (\ref{eq:Objt})-(\ref{eq:FIntg}). To this end, anytime a new global solution is found (i.e., if the search process finds a new solution with an objective function that is lower than the objective function value of the best solution found so far) we resolve model (\ref{eq:Objt})-(\ref{eq:FIntg}) with sprayers routing decision variables being fixed but only small number of refilling points is potentially assessed. \\

More formally, denote the set of ordered refilling points in sprayer $k$ route as $\mathcal{R}_k = \lbrace i_1^k, i_2^k,...\rbrace$. Note that node $i_1^k$ is the first refilling point for sprayer $k$ within its route, node $i_2^k$ is the second refilling point for sprayer $k$ within its route, and so on. For any positive integer $\kappa$, denote $i_{1,-\kappa}^k$ as the $\kappa^{th}$ node \textcolor{red}{preceding} node $i_1^k$ in sprayer $k$ route and denote $i_{1,+\kappa}^k$ as the $\kappa^{th}$ node succeeding node $i_1^k$ in sprayer $k$ route; the same concept can be applied for any refilling point $i_{j,+\kappa}^k$. Lastly, denote the set of refilling points for any feasible solution as $\mathcal{R}=\mathcal{R}_1\cup\mathcal{R}_2\cup...\cup\mathcal{R}_k$, the set of refilling points for any feasible solution can be calculated following the procedure outlined in Section \ref{sec:serviceTime}. Our idea of the intensive local search is that nodes that are in proximity to nodes of set $\mathcal{R}$ are good candidates as refilling points. We define this proximity as being a successor or a predecessor of members of set $\mathcal{R}$. Denote the set of nodes that are candidates to be refilling points as $\mathcal{C}$. Depending on the value of $\kappa$, we can add nodes to set $\mathcal{C}$ to be explored as refilling points while adding a constraint to mathematical model (\ref{eq:Objt})-(\ref{eq:FIntg}) to exclude other nodes from being considered as refilling points. Namely, the following constraints are added:

\begin{equation}
    \delta_i = 0 \quad \forall i \in \mathcal{N}^f/\mathcal{C}
\end{equation}

\begin{figure}[h]
\centering
\includegraphics[width=14cm]{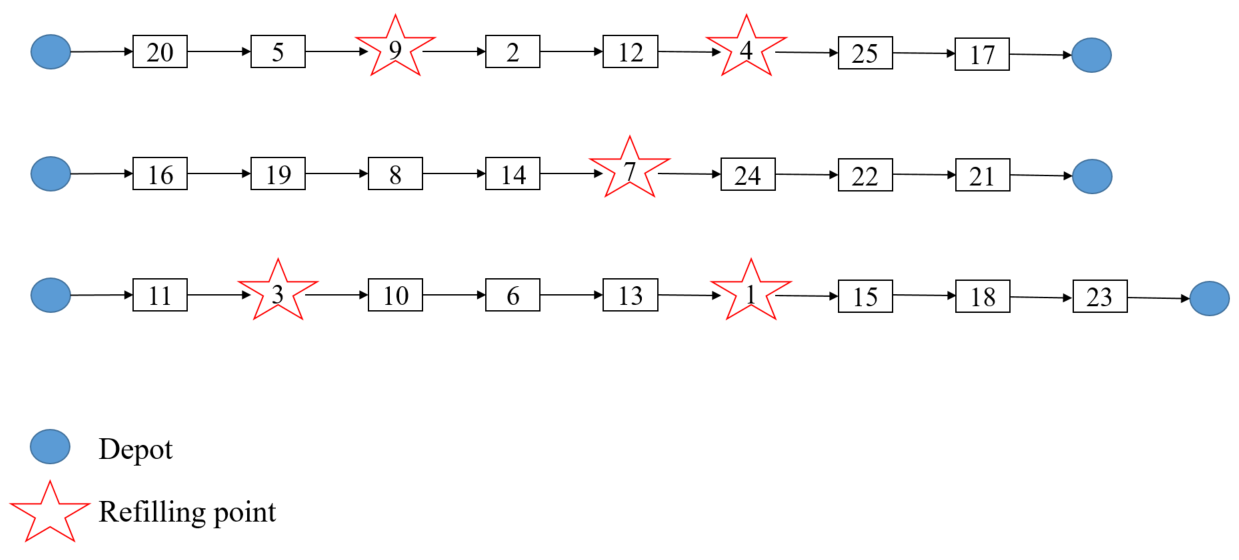}
\caption{Example to illustrate the concept of creating a candidate set of nodes for the local search}
\label{fig:GLS}
\end{figure}

To illustrate the concept of creating set $\mathcal{C}$ we provide a small example based on a feasible solution depicted in Figure \ref{fig:GLS} for an instance of 25 nodes and 3 sprayers. In that feasible solution, sprayer 1 serves nodes (in the following order) 20, 5, 9, 2, 12, 4, 25, and 17, sprayer 2 serves nodes 16, 19, 8, 14, 7, 24, 22, and 21, and sprayer 3 serves nodes 11, 3, 10, 6, 13, 1, 15, 18, and 23. In this feasible solution, the refilling points within the route of sprayer 1 are nodes 9 and 4; hence $\mathcal{R}_1=\lbrace 9,4\rbrace$, for sprayer 2 we have $\mathcal{R}_2=\lbrace 7\rbrace$, and for sprayer 3 we have $\mathcal{R}_3=\lbrace 3,1\rbrace$. If we set the value of $\kappa$ to 1 then the neighborhood of refilling point 9 are nodes 5 and 2 (the nodes proceeding and succeeding node 9 in its position within a route). If we set the value of $\kappa$ to 2 then the neighborhood of node 9 is the set of nodes 20, 5, 2, and 12. We assume that nodes in the neighborhood of the refilling points found using the procedure outlined in Section \ref{sec:serviceTime} present good candidates to be refilling points and hence the process of searching for better refilling points depends on eliminating the possibility of considering nodes other than the ones in the set of  that set $\mathcal{C}$. In this example, if we set the value of $\kappa$ to 0 the set of candidate refilling points is $\mathcal{C}=\lbrace 1,3,4,7,9\rbrace$. On the other hand, if we set the value of $\kappa$ to 1 the set of candidate refilling points is $\mathcal{C}=\lbrace 1,2,5,3,4,7,9,10,11,12,13,14,15,24,25\rbrace$. Clearly, the size of set $\mathcal{C}$ increases as the value of $\kappa$ increases, and thus the computational time needed to perform the intensification procedure. 

\subsubsection{Acceptance Criteria}
\label{sec:acceptCriteria}
We use simulated annealing as an acceptance criterion in our ALNS metaheuristic. The new solution, the current solution, the global best solution, and the function to return the objective function value of any passed solution are defined as $\phi_{n}$, $\phi_{c}$, $\phi^*$ and $f(\,)$, respectively. \textcolor{red}{A \textit{new solution} refers to the solution obtained after performing one ALNS iteration on a \textit{current solution}. More specifically, using a current solution as a starting solution at any iteration, the ALNS performs a destroy operation on the current solution resulting in a partially destroyed solution. Following a destroy operation, the ALNS performs a repair operation on the partially destroyed solution yielding a solution called a \textit{new solution}. Therefore, a new solution is the outcome of any ALNS destroy and repair operations.} \\

A new solution $\phi_{n}$ is always accepted, $f(\phi_{n})<f(\phi^*)$. A worst new solution $\phi_{n}$ can be accepted if $f(\phi_{n}) > f(\phi_{c})$ using a probability function defined as $e^{(f(\phi_{n}) - f(\phi_{c}))/Tem}$ and controlled by parameter $Tem$, where $Tem$ is the temperature. The value of $Tem$ decreases after every iteration according to the formula $Tem = C*Tem$, where $C$ is the cooling rate defined between $(0 < C < 1)$. To compute the initial temperature $Tem_{0}$, we use the formula $Tem_{0} = \rho * f(\phi_{0})$, where $\phi_0$ is the initial feasible solution. \\

The stopping criteria of the ALNS in our work is based on the total number of iterations. 

\subsubsection{Handling infeasible solutions}
\label{sec:infeasibleSol}

A solution generated during the execution of the ALNS might be infeasible due to the violation of constraints (\ref{eq:Ten}) or (\ref{eq:TwentyOne}), namely a solution might violate the requirement that the workload of a sprayer should not exceed $tMax$ or that the waiting time of a sprayer is not zero, respectively. During our design of the ALNS metaheuristic we observed that a big part of the ALNS search process yields infeasible solutions due to violating constraints (\ref{eq:TwentyOne}). Instead of distracting infeasible solutions due to violation of constraints (\ref{eq:TwentyOne}), we allow infeasible solutions temporarily during the search process. On the other hand, violation of constraints (\ref{eq:Ten}) yields infeasible solutions that are not considered during the search process and any solution violating constraints (\ref{eq:Ten}) is discarded. The generation of temporarily infeasible solutions enables a better traversing of the search space because it reduces the chance of getting stuck in a local optimum \cite{cordeau2001unified}, and by alternating between feasible and infeasible regions. Using the appropriate penalty value, we add a penalty term to the objective function to reward solutions that are feasible. To this end, the objective function value that we use during the ALNS execution is:

\begin{eqnarray}
\label{eq:ObjPen}
\displaystyle{\min}&& \displaystyle{\sum_{i\in \mathcal{N}}\sum_{j\in \mathcal{N} i\neq j} t_{ij}(\sum_{k\in\mathcal{K}}x_{ij}^k+ g_{ij}) + \sum_{i\in \mathcal{N} i} \xi \delta_{i}-\sum_{k\in\mathcal{K}}\sum_{i\in\mathcal{N}^f}s_i^k+\lambda\sum_{i\in\mathcal{N}^f}\sum_{k\in\mathcal{K}}m_i^k}
\end{eqnarray}

where the value of $\lambda$ is set to the value of 10. We don't update the value of $\lambda$ dynamically.

\subsubsection{Updating the weights}
\label{sec:weightUp}
In the ALNS metaheuristic, after a predetermined number of iterations the weights of the destroy operators (selection probabilities) update according to the performance of each operator. Because we only design two insert operators, the weight of each insert operator does not change; rather, the weight of the destroy operators $d \in \mathcal{D}$ is updated. \\

The weights of each destroy operator are calculated based on the performance of each operator in during a set of iterations called a \textit{segment}. By using the weights of the previous segment and scores obtained during the current segment the score $\chi_d$ is computed. The formula for updating the value of $\chi_d$ is as follows:
$$\chi_d^{seg-1} = \mu\chi_d^{seg-1}+ (1-\mu)\varphi$$
where $\mu\in(0,1)$ is a smoothing parameter and $\varphi$ is the performance score of a destroy operator and it is assigned as:

\begin{equation*}
  \varphi = \left \{
  \begin{aligned}
    &\varphi_1, && \text{if the new solution is better than the best solution found so far}\\
    &\varphi_2, && \text{if the new solution is better than the current solution}\\
    &\varphi_3, && \text{if the new solution is accepted}\\
    &\varphi_4, && \text{if the new solution is rejected}
  \end{aligned} \right.
\end{equation*}

with $\varphi_1\geq \varphi_2\geq \varphi_3\geq \varphi_4$, in our implementation we use $\varphi_1=7, \varphi_2=4, \varphi_3=2,$ and  $\varphi_4 = 1$. The initial weight $\chi_d$ at the first iteration is set to 1 for all operators. At the end of each segment after updating the values of $\chi$, the probability assigned to each operation $d$ is $w_d=\frac{\chi_d}{\sum_{d\in \mathcal{D}}\chi_d}$. In our work, a segment is set to 200 iterations.

\subsubsection{ALNS Framework: Pseudo code.\label{sec:pseudo code}}

\begin{algorithm}[H]
  \renewcommand{\arraystretch}{0.5}
  \small
\DontPrintSemicolon
  \KwInput{An initial feasible solution $\phi$ produced by the initialization heuristic}
  \KwInput{A set of ALNS parameters: $Tem,\ \mathbf{\varphi},\ MaxIter,\ MaxNoImprov$}
  \KwOutput{The best known solution $\phi^*$}
  Set $iter \gets 0,\ \phi^* \gets \phi,\ \phi_{c} \gets \phi $, $\chi = [1,...,1]$\;
  \While{$iter < MaxIter$}
   {
   		Select a destroy operator $(d)$ and repair operator $(r)$ by a roulette wheel mechanism using weights $\mathbf{w}$\;
     Select a repair operator at random\;
   		Apply destroy operator $d$ on $\phi_c$ and get $\phi^\prime$\;
   		Apply repair operator on $\phi^\prime$ and get $\phi_{n}$\;
   		\If{$f(\phi_{n}) < f(\phi^*)$}
    {
        $NoImproveCounter \gets 0$\;
        $\phi_{n} \gets local-search(\phi_{n})$\ \tcp*{Apply intensive local search to improve the new best solution}
        $\phi^* \gets \phi_{n},\ \phi_{c} \gets \phi_{n}$\;
        \label{line8}
       
    }
    \Else
    {     $NoImproveCounter \gets NoImproveCounter + 1$\;
    	  \If{$f(\phi_{n}) < f(\phi_{c})$}
    {
        $\phi_{c} \gets \phi_{n}$\;
        \label{line13}
        
    }
    \Else
    {
    	$probAccept \gets e^{(f(\phi_{n}) - f(\phi_{c}))/Tem}$\;
    	Generate a random number $U\sim Uniform[0,1]$\;
    	\If{$U < probAccept$}
            {
                $\phi_{c} \gets \phi_{n}$    \tcp*{Accept the new solution which is worst than the current solution for diversification}

            }
            \Else
            {
            	$\phi_c \gets \phi_{c}$ \tcp*{Reject the new solution and keep the current solution for the next iteration}
            	
            }
    }
    }
    $iter\gets iter+1$\;
    $Tem\gets C*Tem$\;
    \If{$iter \% segmentLength == 0$}
        {Update $\mathbf{\chi}$ and calculate $\mathbf{w}$  \;}

    \If{$NoImproveCounter \% MaxNoImprov == 0$}
        {$\phi_c \gets \phi^*,\ NoImproveCounter \gets 0$}
   }
\caption{Overview of ALNS metaheuristic}\label{ALNSalgo}
\end{algorithm}

\subsection{\textcolor{red}{Improvement} Phase}
\label{sec:phase3}
The \textcolor{red}{improvement} phase of our matheuristic focuses on exploring promising solutions based on the search completed in the ALNS phase. Anytime a new best solution is found in phase 2, we save the active sprayers' routing variables (namely, we save all the arcs where $x_{ij}^k=1$ in a pool and we call that pool the high-quality arcs). At the completion of phase 2 and using the high-quality arcs, we resolve model (\ref{eq:Objt})-(\ref{eq:FIntg}) by setting all arcs that do not belong to the set of high-quality pool arcs to 0 for the sprayers only (i.e., $x_{ij}^k=0$ if arc $(i,j)$ is not an element of the high-quality arcs while we leave decision variables $g_{ij}$ unaffected). Our numerical results suggest that as much as 80\% of potential arcs are eliminated from the search. This allows for some nodes to be potentially exchanged as all used parts of the best solutions can be combined, and the model can optimize and merge pieces of different solutions.

\section{Computational Analysis}
\label{sec:Computational}

The aims of the computational analysis in our work are 1) to demonstrate the effectiveness of the matheuristic we propose in our work in terms of finding high-quality solutions in reasonable computational time, more specifically we demonstrate the value of the new elements we introduced in the design and implementation of the matheuristic; 2) to analyze the impact of different variations in modeling on various performance measures, specifically what are the impacts of change in the preference of the decision maker on several performance measures; and 3) to illustrate the value of using an optimization approach considering the integration of spraying and refilling operations as opposed to the approach implemented in practice which solves the problem sequentially. Section \ref{sec:dataGen} provides details of our instances and the parameter setting of the the matheuristic. Section \ref{sec:algoFeatures} reports the detailed results of our proposed matheuristic algorithm features. Section \ref{sec:ALNSVSSolver} demonstrates the performance of our matheuristic against the state-of-the-art solver. An indepth analysis of the effects of using an optimization approach are presented in Section \ref{sec:comparing}.\\

Our algorithms were coded in Python using Jupiter Notebook as an interface and we use Gurobi 9.1.1 optimization software as the MIP solver. All the experiments were conducted on a computer with Intel(R) Core(TM) i7-8700 CPU @ 3.20GHz processor with 16.0 GB RAM.

\subsection{Instance generation}
\label{sec:dataGen}
The problem we study in this paper is a new variation of the VRP with multiple synchronized constraints motivated by a real-world problem faced by one of the largest agriculture companies worldwide. We do have a nondisclosure agreement with this company and therefore we will not be sharing any details regarding the technologies used and the characteristics of the farms they serve. Instead, we will generate data to simulate the problem setting, only in Section \ref{sec:comparing} we present results using real-world data but we do not provide details of that data.\\

Instances are created based on random data. To account for different farm areas and to ensure the reliability of our results, we developed an instance generator to create various problem instances. Four components define an instance: (a) the area of a farm and the number of locations to be sprayed, (b) the number of sprayers, and (c) the service time at each node. The value of various parameters used by our model are as follows:

\begin{itemize}
    \item Number of nodes in the graph: we distinguish three farms, which are defined by the number of nodes to be sprayed and their area. Small farms have typically 15-25 locations (nodes) to be sprayed, medium farms have typically 25-40 nodes, and large farms have more than 40 nodes to be sprayed. The area of a small farm is usually 500 acres, a medium farm is 600-1500 acres, and large farms are 2,000 acres. \textcolor{red}{We set the value of parameter $rd$ used in destroy operator 6 to be 4.00, 4.90, and 5.60 for small, medium, and large-size farms, respectively.}
    \item Parameters $qMin_i$ are set to be random numbers generated in the interval of [1.5-3.5], $qMax_i$ are set to be 2.5 times $qMin_i$ for each node.
    \item Parameter $\xi$ is assumed to be 3 units of time.
    \item Parameter $Q^s$ is assumed to be 15 units of fertilizer.
    \item Parameter $\eta^{sp}$ is assumed to be 2 units of time per unit of fertilizer.
    \item Parameter $\beta$ is assumed to be 2.
\end{itemize}

In practice, for small and medium farms two or three sprayers will be made available while three or four sprayers are usually made available for large farms.\\

Matheuristic parameter settings are as follows:
\begin{itemize}
    \item We set the maximum number of iterations for the ALNS metaheuristic to be 200*$\vert \mathcal{N}^f\vert$.
    \item The time limit to solve the mathematical model in phase 3 is 2 hours for small instances, 4 hours for medium instances, and 6 hours for large instances. Our earlier computational results show that solving the mathematical model of phase 3 takes time for large instances and setting up the time limit of the solver to 2 or 4 hours is not enough.
    \item The time limit to solve the mathematical model of the intensive local search is set to 60 seconds. While this does not always yield an optimal solution with a 0.0\% optimality gap, our preliminary results suggest that 60 seconds provides a good trade-off between getting a high-quality solution and computational time. 
    \item The number of iterations before updating the weights of the destroy operators is set to be 200 iterations. 
    \item The initial temperature is set to be 100 times the objective function value of the initial solution and the cooling rate is set to 0.995. 
\end{itemize}

\subsection{Evaluation of the new features}
\label{sec:algoFeatures}
During the design of our matheuristic, we introduced new features such as the way of performing the construction phase to find an initial solution, new destroy operators, and intensifying local search during the execution of the ALNS phase. This section presents a comprehensive analysis to evaluate the impact of each of these new features on the performance of the matheuristic. 

\subsubsection{Construction phase}
In section \ref{sec:feasible solution}, we introduced two techniques to find an initial solution (i.e., phase 1 of our matheuristic) to be used by the ALNS metaheuristic as a starting point. To assess the quality of each technique, we compare the solution quality after running phase 1 and phase 2 under each technique. We only use medium instances as a way to assess the performance of the construction algorithm.\\

The results of the comparison are displayed in Table \ref{tab:initialization}. The first column reports the number of nodes (locations to be sprayed), the second column shows the number of sprayers, the third set of columns reports the gap between the lower bound and the best value reported at the termination of the ALNS and the computational time in seconds under the greedy insert, and similarly the fourth set of columns reports on the gap between the lower bound and the best value reported at the termination of the ALNS and the computational time in seconds under the technique of the cluster-first-route-second. The gap is calculated as $(\frac{z^*-LB}{LB})*100\%$ where $z^*$ is the best value found by the ALNS metaheuristic and $LB$ is the lower bound value calculated following the procedure outlines in Section \ref{sec:LowerBound}. Values in bold font represent the best performance. Each row in Table \ref{tab:initialization} represents the mean value of a performance measure taken by running 15 simulation replicates under each instance size. We perform 15 runs per instance to capture variations due to input parameters. The only variations between runs of the same instance are the minimum quantity of fertilizer to be applied at each node and the coordinates of nodes in the network.\\

\begin{table}[h]
\begin{center}
\caption{Averaged computational time and gap for the initialization techniques.}
\label{tab:initialization}
\begin{tabular}{llllll}
\multicolumn{1}{c}{$\vert \mathcal{N}^f\vert$} & \multicolumn{1}{c}{$\vert \mathcal{K}\vert$} & \multicolumn{2}{c}{Greedy}                                     & \multicolumn{2}{c}{Clustering}                                 \\
                                             &                                              & \multicolumn{1}{c}{Gap (\%)} & \multicolumn{1}{c}{Time (sec.)} & \multicolumn{1}{c}{Gap (\%)} & \multicolumn{1}{c}{Time (sec.)} \\ \hline
25	&	2	&	\textbf{107.43}	&	 1,966 	&	132.45	&	 1,991 	\\
25	&	3	&	\textbf{123.93}	&	 1,805 	&	173.88	&	 1,838 	\\
25	&	4	&	\textbf{98.30}	&	 1,433 	&	222.88	&	 1,425 	\\
30	&	2	&	74.33	&	 5,386 	&	\textbf{61.96}	&	 5,179 	\\
30	&	3	&	63.30	&	 4,750 	&	\textbf{60.95}	&	 4,665 	\\
30	&	4	&	\textbf{87.84}	&	 3,780 	&	164.78	&	 3,815 	\\
35	&	2	&	\textbf{142.92}	&	 9,764 	&	158.18	&	 9,213 	\\
35	&	3	&	\textbf{111.53}	&	 8,972 	&	173.50	&	 9,039 	\\
35	&	4	&	\textbf{157.98}	&	 7,849 	&	293.89	&	 8,000 	\\
40	&	2	&	\textbf{136.51}	&	 19,936 	&	155.67	&	 18,037 	\\
40	&	3	&	\textbf{126.51}	&	 16,286 	&	178.06	&	 16,875 	\\
40	&	4	&	\textbf{115.35}	&	 14,397 	&	221.27	&	 14,554 	\\ \hline
       & Avg. & 112.16  &  & 166.79 &                          \\
       & Min  & 63.30   &  & 60.95  &                         \\
       & Max  & 157.98  &  & 293.89 &                         \\ \hline
\end{tabular}
\end{center}
\end{table}

From Table \ref{tab:initialization} we note that the greedy initialization technique demonstrates better performance than the clustering technique. On the other hand, we observe that the CPU time of both techniques has similar values. Therefore, we conclude that the greedy insert technique is a better technique to generate the initial feasible solution for our matheuristic, hence we will be using the greedy insert technique as a means of generating the initial feasible solution.\\

\begin{figure}[h]
\centering
\includegraphics[width=12cm]{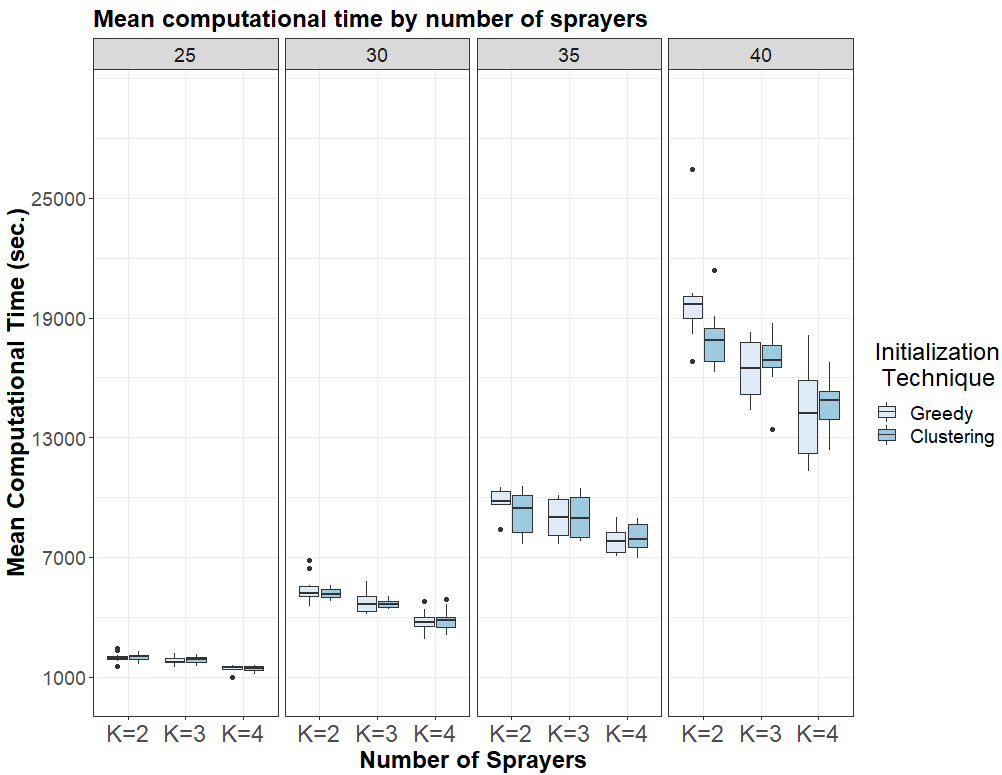}
\caption{Computational time in seconds for each initialization technique under different number of nodes and number of sprayers}
\label{fig:InitTime}
\end{figure}

Lastly, to further illustrate more details on the computational time in seconds under each technique, Figure \ref{fig:InitTime} displays box plots of the computational time in seconds for each initialization technique under different number of nodes and number of sprayers. We observe that the variability in computational time increases as the number of nodes and sprayers increases.

\subsubsection{Value of ALNS with intensive local search}
In section \ref{sec:localSearch}, we introduced a way to find higher quality solutions by implementing a form of local search to find the optimal service time, refilling locations, refilling time schedule, and tanker's routing given fixed sprayers routes using mathematical modeling. In this subsection, we explore the value of utilizing such an approach as a means for intensive local search.\\

Recall that by tuning the value of $\kappa$, the size of the neighborhood of a solution might be explored changes in size. In this section, we demonstrate the value of integrating the intensive local search method we developed in Section \ref{sec:localSearch} under different values of $\kappa$. Table \ref{tab:glsTable} summarizes the outcome of this numerical testing. Each row in Table \ref{tab:glsTable} represents the mean value of a performance measure taken by running 15 simulation replicates under each instance size. We perform 15 runs per instance to capture variations due to input parameters. The only variations between runs of the same instance are the minimum quantity of fertilizer to be applied at each node and the coordinates of nodes in the network. \\

The first column in Table \ref{tab:glsTable} shows the number of nodes and the second column shows the number of sprayers. The third set of columns shows the performance of running the numerical experiments without using the intensive local search, the second column is labeled as \textit{No LS}. The first sub-column shows the computational time in seconds of running Phases 1 and 2 of our matheuristic and the second sub-column shows the gap between the best solution found at the termination of the ALNS phase against the lower bound. Note that we do not use phase 3 in this comparison to highlight the effects of using a local search on the performance of the ALNS phase. The fourth set of columns, labeled as \textit{$\kappa = 0$}, in Table \ref{tab:glsTable} shows the computational time in seconds and the gap in a strategy where the value of $\kappa$ is set to 0, namely the refilling points are established following the procedure in Section \ref{sec:serviceTime}. The fifth set of columns, labeled as \textit{$\kappa=1$}, shows the computational time in seconds and the gap in a strategy where the value of $\kappa$ is set to 1. And finally, the sixth set of columns, labeled as \textit{Hybrid}, shows the computational time in seconds and the gap in a strategy where no intensive local search is implemented during the first one-third of the ALNS iterations, then the value of $\kappa$ is set to 0 for the second one-third of the iterations, and the value of $\kappa$ is set to 1 for the last one-third of the iterations. \\

\begin{table}[h]
\caption{Averaged computational time and gap under each intensive local search strategy.}
\label{tab:glsTable}
\centering
\begin{tabular}{rrcccccccc}
 &  & \multicolumn{2}{c}{No search} & \multicolumn{2}{c}{$\kappa=0$} & \multicolumn{2}{c}{$\kappa=1$} & \multicolumn{2}{c}{Hybrid} \\
$\vert \mathcal{N}^f\vert$ & $\vert \mathcal{K}\vert$ & \begin{tabular}[c]{@{}c@{}}Time \\ (sec.)\end{tabular} & \begin{tabular}[c]{@{}c@{}}Gap\\  (\%)\end{tabular} & \begin{tabular}[c]{@{}c@{}}Time \\ (sec.)\end{tabular} & \begin{tabular}[c]{@{}c@{}}Gap\\  (\%)\end{tabular} & \begin{tabular}[c]{@{}c@{}}Time \\ (sec.)\end{tabular} & \begin{tabular}[c]{@{}c@{}}Gap\\  (\%)\end{tabular} & \begin{tabular}[c]{@{}c@{}}Time \\ (sec.)\end{tabular} & \begin{tabular}[c]{@{}c@{}}Gap\\  (\%)\end{tabular} \\ \hline
25	&	2	&	2,098	&	195.06	&	2,072	&	\textbf{177.12}	&	2,063	&	186.46	&	\textbf{1,998}	&	178.47	\\
25	&	3	&	1,845	&	208.66	&	1,883	&	144.78	&	\textbf{1,821}	&	137.41	&	1,851	&	\textbf{137.29}	\\
25	&	4	&	\textbf{1,322}	&	152.99	&	1,342	&	\textbf{119.35}	&	1,359	&	123.34	&	1,336	&	127.73	\\
30	&	2	&	\textbf{4,428}	&	207.55	&	4,660	&	190.58	&	4,768	&	186.5	&	4,541	&	\textbf{184.76}	\\
30	&	3	&	\textbf{4,131}	&	241.34	&	4,359	&	167.33	&	4,188	&	\textbf{164.89}	&	4,141	&	185.11	\\
30	&	4	&	\textbf{3,540}	&	281.88	&	3,632	&	154.9	&	3,694	&	\textbf{148.72}	&	3,557	&	151.24	\\
35	&	2	&	\textbf{8,704}	&	219.75	&	9,393	&	205.73	&	9,086	&	\textbf{191.87}	&	8,915	&	198.58	\\
35	&	3	&	8,478	&	206.49	&	8,801	&	173.21	&	8,443	&	\textbf{166.75}&	\textbf{8,128}	&	173.14	\\
35	&	4	&	6,964	&	216.29	&	7,092	&	\textbf{189.15}	&	7,234	&	190.73	&	\textbf{6,795}	&	191.23	\\
40	&	2	&	\textbf{18,394}	&	244.55	&	18,515	&	\textbf{200.84}	&	19,548	&	203.13	&	19,459	&	216.35	\\
40	&	3	&	\textbf{15,278}	&	194.68	&	15,730	&	\textbf{181.18}	&	16,363	&	181.43	&	16,115	&	192.49	\\
40	&	4	&	14,792	&	180.86	&	13,978	&	\textbf{154.1}	&	14,444	&	156.55	&	\textbf{14,272}	&	158.33	\\ \hline
\multicolumn{1}{l}{} & Avg. &  & 212.51 &  & 171.52  &  & \textbf{169.82} &  & 174.56  \\
\multicolumn{1}{l}{} & Max. &  & 281.88 &  & 205.73 &  & \textbf{203.13} &  & 216.35  \\
\multicolumn{1}{l}{} & Min. &  & 152.99 &  & \textbf{119.35} &  & 123.34 &  & 127.73 \\ \hline
\end{tabular}
\end{table}

Using the values provided in Table \ref{tab:glsTable} we observe that, on average, the computational time does not change significantly when using different intensive local search strategies. On the other hand, we observe that the solution quality, on average, improves as the value of $\kappa$ increases. This is indeed consistent with what anyone would expect, as the size of the neighborhood explored increases, more potential should be realized to find a higher-quality solution. More specifically, the average improvement in the computational gap as calculated against the lower bound using the local search where the value of $\kappa$ is set to 1 when compared to no search strategy is $6.45\%$. On the other hand, the average improvement in the computational gap as calculated against the lower bound using the local search where the value of $\kappa$ is set to 1 when compared to local search with $\kappa=0$ strategy is $3.30\%$. \textcolor{red}{From the results shown in Table \ref{tab:glsTable} we observe that the computational times of the intensive local search are not substantially longer than the computational times without the local search. This can be explained by pointing out that the computational time required to solve the MILP model with fixed routing variables for sprayers and fixing a subset of the nodes as candidate locations for refilling is solved in less than 5 seconds for the largest instances and less than 2 seconds for medium instances. Such short computational time needed to perform the local search does not substantially increase the computational times. We also observed during the computational analysis that the number of times a new best solution is found using any of the intensive local search strategies is less than the number of times a new best solution is found under the `no search' strategy.} \\

\begin{figure}[!h]
\centering
\includegraphics[width=14cm]{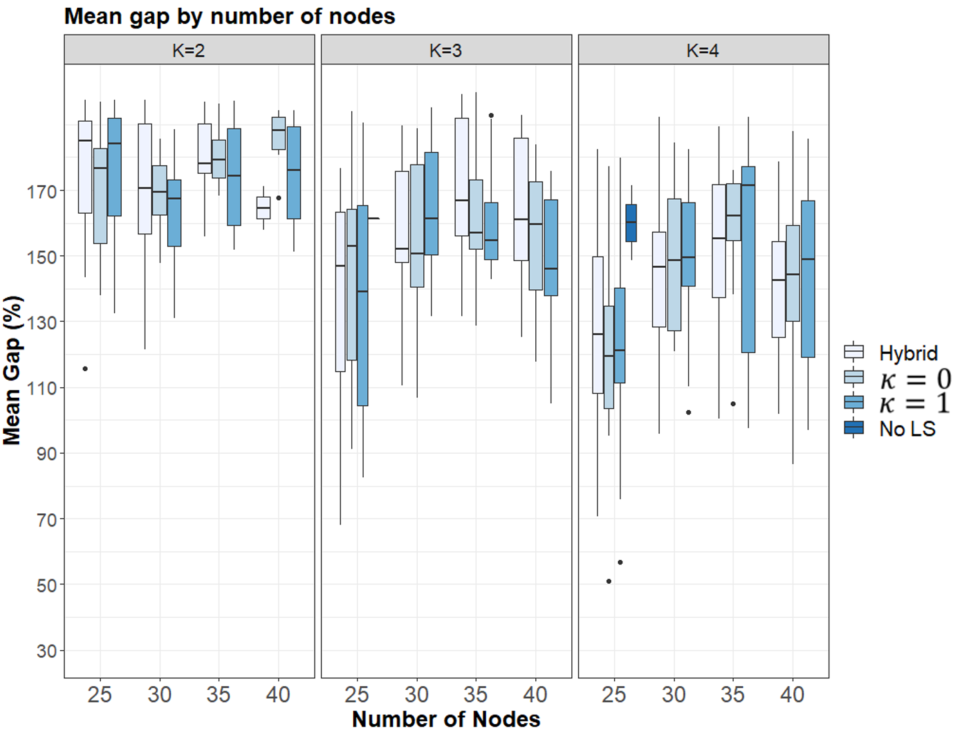}
\caption{Computational gap for each local search strategy}
\label{fig:glsGap}
\end{figure}

\begin{figure}[!h]
\centering
\includegraphics[width=14cm]{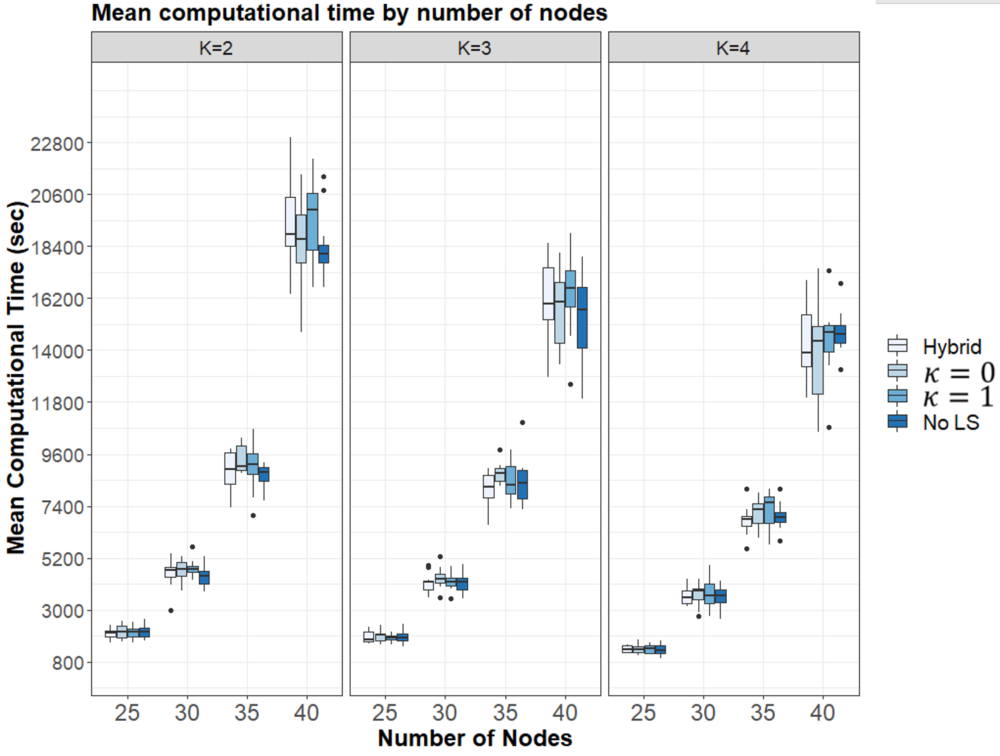}
\caption{Computational time for each local search strategy}
\label{fig:glsTime}
\end{figure}

Lastly, Figures \ref{fig:glsGap} and \ref{fig:glsTime} illustrate the mean gap and computational times under each strategy, respectively. We provide Figures \ref{fig:glsGap} and \ref{fig:glsTime} to show the variation in the results as Table \ref{tab:glsTable} only reports the average. \textcolor{red}{Using the boxplot displayed in Figure \ref{fig:glsGap} we observe that the gap between phase 2 and lower bound is always larger than the gap between phase 3 and the lower bound. On the other hand, from the boxplot displayed in Figure \ref{fig:glsTime} we observe that the computational time increases exponentially as the number of nodes increases but never exceeds the 10 hours time limit used in practice.} \\

\subsubsection{Value of new destroy operators}
In this subsection, we show the effectiveness of the destroy operators designed on the performance of the ALNS. Recall that the destroy operators are selected probabilistically. Operators with higher performance have a better probability of being selected. This implies that the more efficient an operator has been in prior iterations, the more likely it will be chosen in subsequent iterations. Thus; we assess the effectiveness of a destroy operator by examining the probability of selecting this operator over iterations. For the purposes of this analysis, we show the results on instances with 35 nodes and 4 sprayers, it is worth mentioning that we observe similar behavior across different problem instances. We perform 15 independent runs and we provide the average values.\\

\begin{figure}[!h]
\centering
\includegraphics[width=14cm]{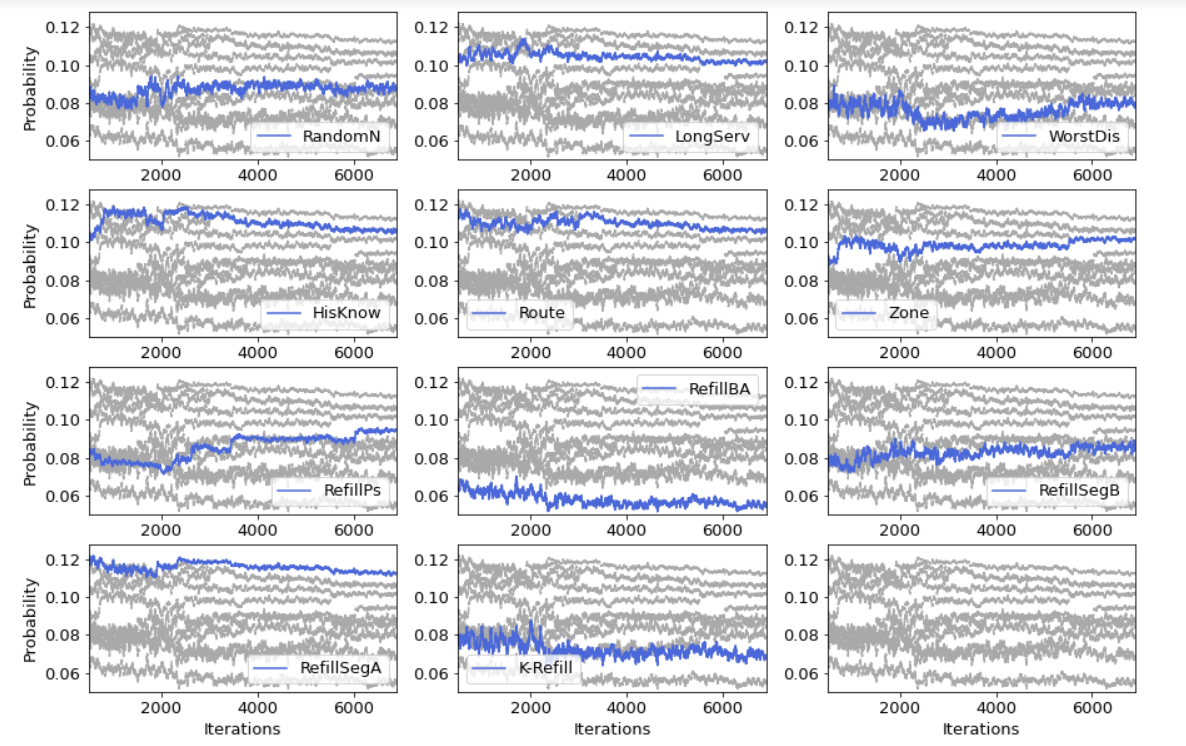}
\caption{Probability of destroy operator selection over a number of iterations}
\label{fig:ValueOp}
\end{figure}

Figure \ref{fig:ValueOp} shows the average probability of a destroy operator being chosen over 7,000 iterations (recall that the number of iterations of the ALNS is set to 200*$\vert \mathcal{N}\vert$). The last sub-plot in Figure \ref{fig:ValueOp} on the lower right shows the performance of all destroy operators without any highlight.\\

In each subplot, we present the probability of an operator's selection averaged across all runs for the selected instance size. The operators we use from the literature are labeled as: random - $RandomN$, route - $Route$, worst distance - $WorstDis$, historical knowledge - $HisKnow$, and zone - $Zone$. The destroy operators we designed for this specific study are labelled as: refill Position - $RefillPs$, neighbors of refill positions - $RefillBA$, longest distance - $LongServ$, subtour prior to refilling position - $RefillSegB$, subtour following a refilling position - $RefillSegA$ and K-Neighbors of refill positions - $K-Refill$. See Section \ref{sec:Destroy} for a detailed description of each of the destroy operators. \textcolor{red}{We also report information on the performance of each operator in contributing to finding a new best solution by reporting the ratio of the number of times a new best solution was found by each operator divided by the total number of times a new best solution is found. Following the order of the destroy operators provided in this paragraph, these values are: 0.17, 0.14, 0.03, 0.10, 0.09, 0.20, 0.06, 0.03, 0.03, 0.15, and 0.02, respectively. We observe, on average, that the zone removal operator provides the largest contribution to finding a new best solution, followed by random removal, subtour following a refilling position, and historical knowledge.}\\

As seen in Figure \ref{fig:ValueOp}, the majority of the operators we use from the literature perform well. Overall three of the destroy operators we introduce, $RefillSegA$, $RefillSegB$ and $RefillPs$ made significant contributions while the remaining two did not perform well. We observe that both subtour destroy operators we use performed well with $RefillSegA$ showing the highest probability of selection out of all the destroy operators we use. In terms of ranking the probability of selecting an operator, from Figure \ref{fig:ValueOp} we can see that historical knowledge, subtour following a refill, and route destroy operators have the highest probability of being selected, the next operators are longest service, zone, and refill position destroy operators have a high probability of being selected. Lastly, we observe that k-refill and neighbors of refill positions destroy operators perform poorly. \\

\subsection{Results of comparing the performance of the matheuristic against the state-of-the-art solver}
\label{sec:ALNSVSSolver}

We now present the average results obtained by solving mathematical model (\ref{eq:Objt})-(\ref{eq:FIntg}) presented in Section \ref{sec:model} using Gurobi solver and compare it with the results obtained by using our developed matheuristic. To further highlight the value of having phase 3 in our developed matheuristic, we report detailed results on the computational time and solution quality at the end of phase 2 and phase 3. Table \ref{tab:solverVSAlgo} presents the results for instances with a number of nodes 15, 17,$\dots$, 25 and for a number of sprayers equals 2 and 3. \\

The first column of Table \ref{tab:solverVSAlgo}, labeled as \textit{Ins.}, we provide information about the instance, namely the number of nodes and the number of sprayers in the form of $\vert\mathcal{N}^f\vert-\vert\mathcal{K}\vert$. The second column, labeled as \textit{Phase 2}, reports information regarding phase 2 of the matheuristic (i.e., at the termination of the ALNS metaheuristic). The first sub-column reports the mean computational time in seconds, the second sub-column displays the mean objective function value (denoted as $z_{Ph2}$), and the mean gap between the objective function value and the objective function value of the lower bound is reported in the third sub-column. More specifically, $Gap=\frac{z_{Ph2}-z^*}{z^*}*100\%$, where $z^*$ is the objective function value calculated using the lower bound. Then the next subset of columns, labeled \textit{Phase 3}, reports the mean computational time in seconds of phase 2, the mean objective function value (denoted as $z_{Ph3}$), and the mean gap between the objective function value and the objective function value of the lower bound (i.e., $Gap=\frac{z_{Ph3}-z^*}{z^*}*100\%$). The next set of columns, labeled as \textit{Solver}, displays the mean objective function value as found by Gurobi solver at the termination, the gap as provided by the solver, and the gap between the objective function value provided by the solver and the lower bound. Note that we do not report the computational time of the solver since it always runs out of time. The time limit we set for Gurobi is 4 hours. The fifth column shows the mean gap between the solver and the matheuristic calculated as $Gap=\frac{z_{Sol}-z_{Ph3}}{z_{Sol}}*100\%$. Lastly, the sixth column shows the value of the objective function under each instance using the lower bound procedure outlined in Section \ref{sec:LowerBound}.\\

Results in Table \ref{tab:solverVSAlgo} show how the developed matheuristic provides a higher quality solution compared to solving the same problem using Gurobi solver in a fraction of time compared to the solver as well. \\

\begin{table}[h]
\caption{Average results of the computational results of the matheuristic, solver, and the lower bound.}
\label{tab:solverVSAlgo}
\begin{tabular}{llllllllllcl}
 & \multicolumn{3}{c}{Phase 2} & \multicolumn{3}{c}{Phase 3} & \multicolumn{3}{c}{Solver} & Sol. vs Mathe. & \multicolumn{1}{c}{LB} \\ \cline{2-10}
Ins. & \multicolumn{1}{c}{\begin{tabular}[c]{@{}c@{}}Time\\  (sec.)\end{tabular}} & \multicolumn{1}{c}{$z_{Ph2}$} & \multicolumn{1}{c}{\begin{tabular}[c]{@{}c@{}}Gap \\ (\%)\end{tabular}} & \multicolumn{1}{c}{\begin{tabular}[c]{@{}c@{}}Time\\  (sec.)\end{tabular}} & \multicolumn{1}{c}{$z_{Ph3}$} & \multicolumn{1}{c}{\begin{tabular}[c]{@{}c@{}}Gap \\ (\%)\end{tabular}} & \multicolumn{1}{c}{$z_{Sol}$} & \multicolumn{1}{c}{\begin{tabular}[c]{@{}c@{}}Solv.\\  Gap(\%)\end{tabular}} & \multicolumn{1}{c}{\begin{tabular}[c]{@{}c@{}}Gap \\ (\%)\end{tabular}} & Gap (\%) & \multicolumn{1}{c}{$z^*$} \\ \hline
15-2&204&24.22&162.05&	1.56&24.22	&	162.03	&	22.01	&	86.11	&	134.88	&	-11.30	&	9.50	\\
17-2&297&27.25&149.09&	4.26&27.06	&	146.96	&	25.35	&	97.37	&	131.36	&	-6.80	&	11.09	\\
19-2&538&30.52&152.24&	6.48&30.45	&	151.79	&	29.20	&	104.04	&	142.26	&	-3.99	&	12.06	\\
21-2&991&28.92&149.03&	32.88&28.86	&	148.47	&	28.73	&	119.98	&	145.99	&	-2.35	&	11.64	\\
23-2&1459&31.48&148.60&	50.82&31.27	&	146.92	&	32.83	&	133.54	&	158.99	&	3.58	&	12.73	\\
25-2&1944&36.02&164.98&	202.98&35.82	&	163.39	&	37.07	&	136.42	&	172.33	&	2.34	&	13.64	\\
15-3&145&21.45&100.36&	1.56&21.08&	96.45	&	19.46	&	53.80	&	79.14	&	-9.46	&	10.89	\\
17-3&238&18.71&87.72&	1.95&18.65&	87.18	&	16.99	&	99.70	&	70.25	&	-10.10	&	9.94	\\
19-3&446&19.18&93.80&	3.30&19.14&	93.37	&	19.50	&	124.03	&	97.24	&	1.60	&	9.84	\\
21-3&803&28.31&124.03&	29.63&28.20&	123.21	&	29.36	&	114.75	&	131.71	&	3.08	&	12.63	\\
23-3&1229&30.10&144.03&	174.04&30.02&	143.47	&	35.14	&	125.71	&	185.25	&	14.12	&	12.31	\\
25-3&1843&32.86&155.24&	49.73&32.81&	154.75	&	37.17	&	136.82	&	186.35	&	9.36	&	12.94	\\ \hline
Avg.	&	845	&		&	135.93	&	46.60	&		&	134.83	&		&	111.02	&	136.31	&	-0.83	&		\\ 
Max.	&	1944&		&	164.98	&	202.98	&		&	163.39	&		&	136.82	&	186.35	&	14.12	&		\\
Min.	&	145	&		&	87.72	&	1.56	&		&	87.18	&		&	53.80	&	70.25	&	-11.30	&		\\
 \hline
\end{tabular}
\end{table}

Several interesting observations can be drawn from the analysis of the results from Table \ref{tab:solverVSAlgo}. First, we observe that compared to the MIP solver, our matheuristic can in almost all cases find better solutions. The improvements are remarkable as the number of nodes grows as well as the number of sprayers increases. Except for the case with 15 nodes and 3 sprayers in which the average solution obtained worsens slightly (less than 1.50\%). Compared to the solutions obtained by phase 2, we can observe that our ALNS combined with the construction phase can find a higher quality solutions compared to the solver when the number of nodes is larger than 19 for both the cases of 2 or 3 sprayers. This observation can be made by comparing the values of the objective function obtained at the end phase 2 against the objective function value reported by the solver (see the second sub-column of the second column against the first sub-column of the fourth column). \\

Considering the time needed to achieve these improvements, the effectiveness of the proposed matheuristic is even better highlighted. The results from Table \ref{tab:solverVSAlgo} show that the matheuristic approach can provide better solutions considerably faster. It is worth mentioning that with 25 nodes and 3 sprayers, the solver was not able to find a feasible point in 4 instances within the time limits. On the other hand, the average running time of our matheuristic is 2,948 seconds for phases 1 and 2 and 1,260 seconds for phase 3, against more than 14,400 seconds for the MIP solver, which represents a reduction of almost 70\%. These results highlight the importance of the matheuristic approach in both aspects: to find high-quality solutions significantly faster.\\

Even with the great results of our matheuristic, we still observe that the gap between the lower bound and the objective function of the developed matheuristic is significant, specifically in the range of 46-86\%. This observation highlights the fact that the lower bound utilized in our numerical results is rather loose and far from the optimal objective function value of mathematical model (\ref{eq:Objt})-(\ref{eq:FIntg}) due to the fact that the synchronization constraints are relaxed. Lastly, results from Table \ref{tab:solverVSAlgo} demonstrate that adding phase 3 brings a significant improvement in our matheuristic by looking at the significant improvement in the solution quality as reported by Phase 3. This hypothesis can be examined by looking at the computational gap of the MIP solver against the computational gap against the lower bound. For example, in instances with 23 nodes and 3 sprayers, we observe that the solver gap is in fact smaller than the gap against the lower bound.\\

\subsubsection{Analyzing the quality of the lower bound}
\begin{figure}[h]
     \centering
     \begin{subfigure}{0.45\textwidth}
         \centering
         \includegraphics[width=\textwidth]{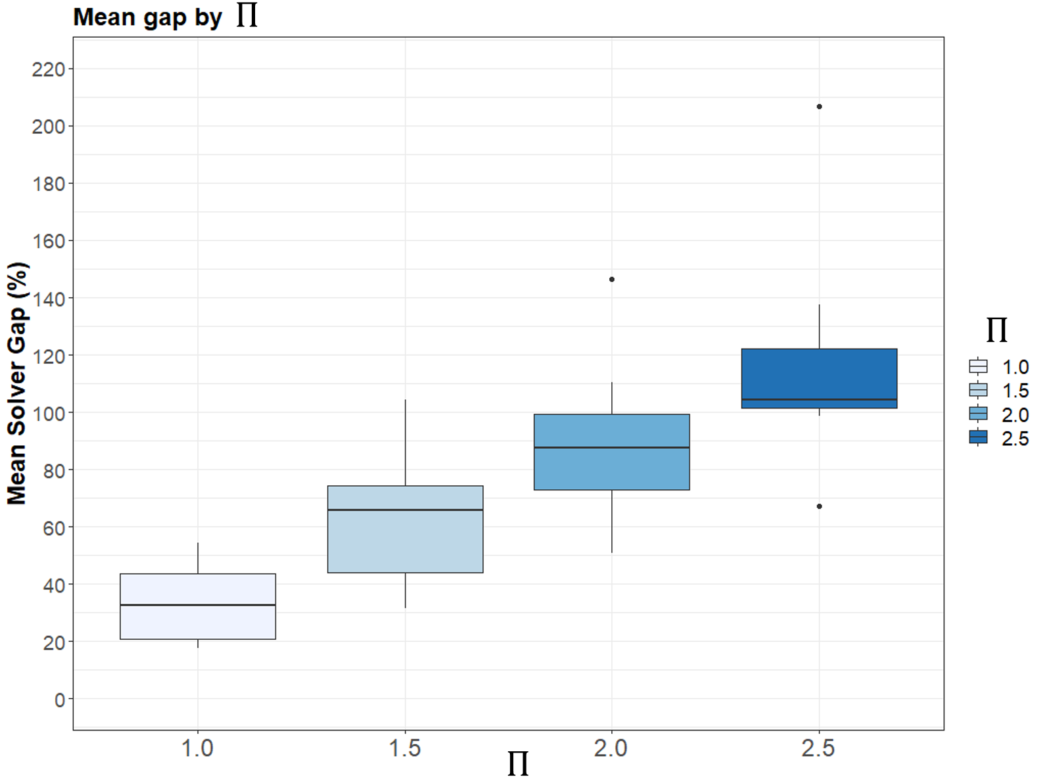}
         \caption{Solver optimality gap}
         \label{fig:solverGap}
     \end{subfigure}
     \hfill
     \begin{subfigure}{0.45\textwidth}
         \centering
         \includegraphics[width=\textwidth]{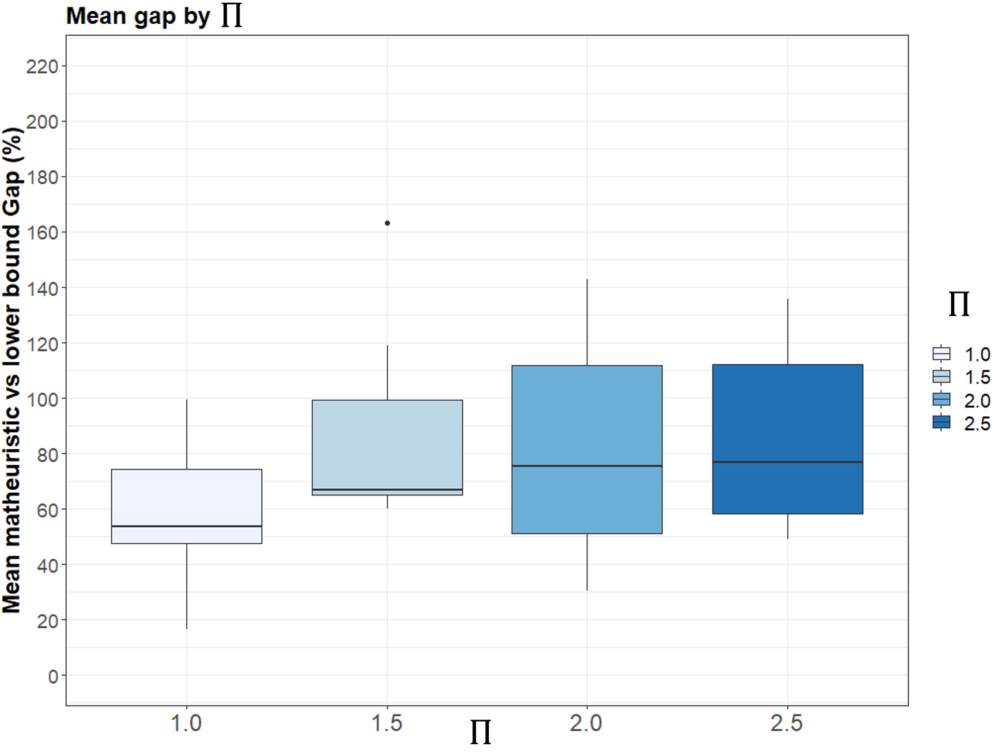}
         \caption{Gap between matheuristic and lower bound}
         \label{fig:mathLb}
     \end{subfigure}
     \vfill
     \begin{subfigure}{0.45\textwidth}
         \centering
         \includegraphics[width=\textwidth]{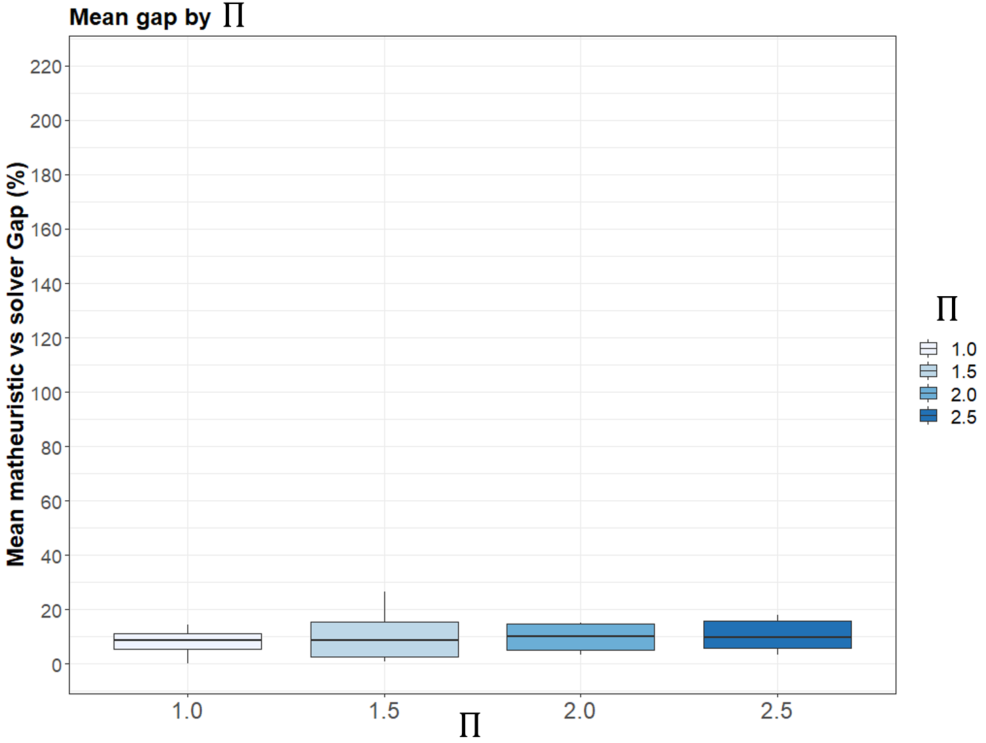}
         \caption{Gap between matheuristic and solver}
         \label{fig:solverMath}
     \end{subfigure}
     \hfill
     \begin{subfigure}{0.45\textwidth}
         \centering
         \includegraphics[width=\textwidth]{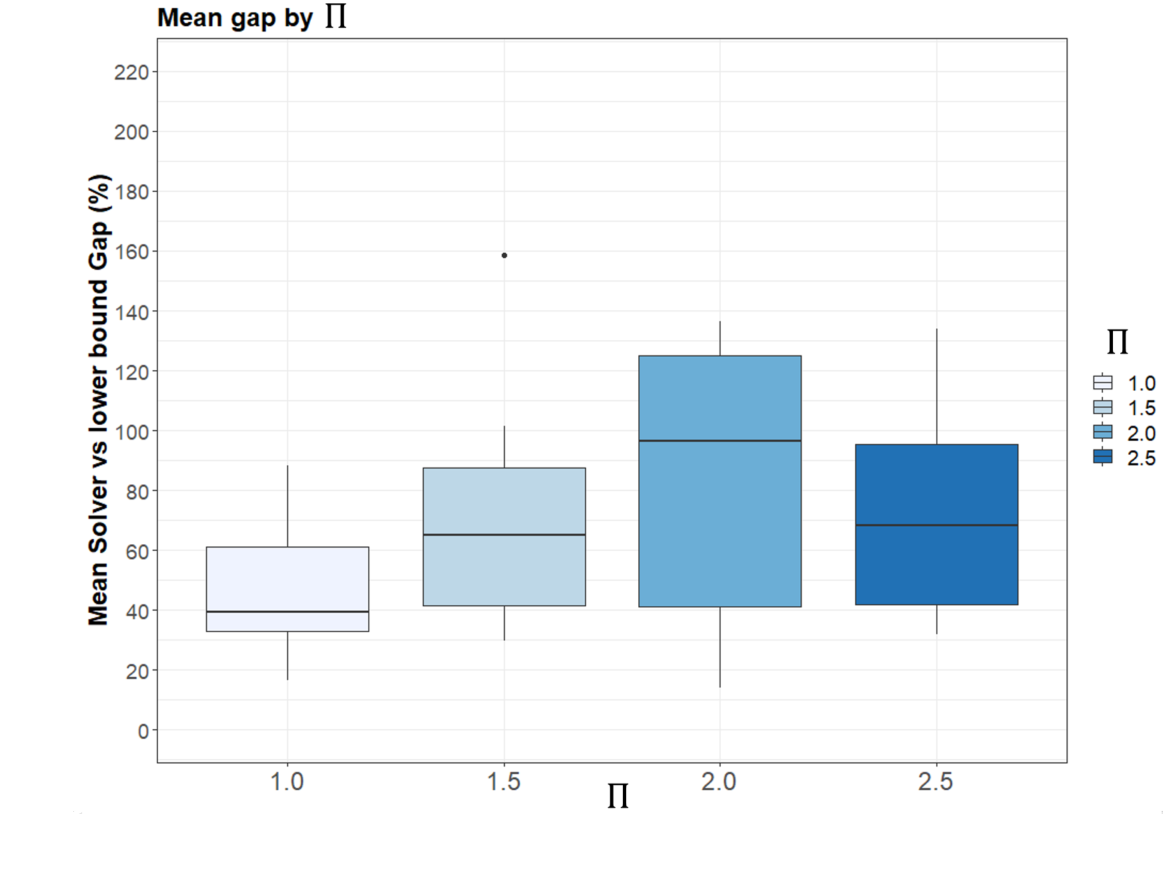}
         \caption{Gap between solver and lower bound}
         \label{fig:gapSolLB}
     \end{subfigure}
        \caption{Solver optimality gap and gap between solver, matheuristic, and lower bound under different values of $\Pi$}
        \label{fig:gaps}
\end{figure}

\textcolor{red}{From the results summarized in Tables (\ref{tab:initialization})-(\ref{tab:solverVSAlgo}) we observe that the values of the optimality gaps reported by the solver at termination and the gap between the lower bound and the best solution found by the developed matheuristic are unusually large, sometimes as large as 160\%. Given the complexity of the problem we study in this work and the fact that the number of binary decision variables grows exponentially with the number of nodes in the network and the number of sprayers, such values of gaps should not be surprising. Furthermore, the total service time positively affects the objective function value by reducing it (see the last term in the objective function (\ref{eq:Objt})). By considering the fact that the range of service time at each node is wide, the value of $qMax_i$ is set to be $2.5*qMin_i$ as discussed in Section \ref{sec:dataGen}. We believe that the wide range of $[qMin_i,qMax_i]$ negatively affects the quality of the lower bound. To test the hypothesis that the low quality of the lower bound and the high values of optimality gaps are partly due to the wide range of service levels at each node, we conduct additional experiments using different values of $\Pi*qMax_i$ where $\Pi\in\lbrace 1.0,1.5,2.0,2.5\rbrace$ and we maintain the fact that the values of $qMin_i$ are randomly generated in the interval of [1.5-3.5]. In summary, we hypothesis that as the value of $qMax_i$ approaches the value of $qMin_i$, the value of the optimality gap reported by Gurobi solver and the quality of the gap between the lower bound and the best solution found by the matheuristic shall improve substantially. We use instances with $\mathcal{N}^f=17$ and $\mathcal{K}=3$ and we generate 10 instances with different locations and $qMin_i$ values. We record the optimality gap of Gurobi solver at termination with a time limit of 4 hours, the value of the lower bound, and the objective function value of the best solution found by the matheuristic. Under different values of $\Pi$ we report the gap values in Figure \ref{fig:gaps}. Figure \ref{fig:solverGap} displays the optimality gap of Gurobi solver, Figure \ref{fig:mathLb} displays the gap between the matheuristic and the lower bound calculated as $\frac{z_{Ph3}-z^*}{z^*}*100\%$, Figure \ref{fig:solverMath} displays the gap between the solver and the best value of the matheuristic $\frac{z_{Ph3}-z_{Sol}}{z_{Sol}}*100\%$, and finally Figure \ref{fig:gapSolLB} displays the gap between the solver and the lower bound calculated as $\frac{z_{Sol}-z^*}{z^*}*100\%$. Clearly, we observe that as the range of $[qMin_i,qMax_i]$ decreases, the value of the optimality gap and the quality of the lower bound improves significantly highlighting the fact that the low quality of the lower bound might be attributed mostly due to the wide range of $[qMin_i,qMax_i]$. Furthermore, we notice that setting the value of $\Pi$ to 1 yields the tightest bounds which demonstrates that the vehicle routing problem with synchronization and variable service time is a much harder problem that the vehicle routing problem with synchronization and fixed service time counterpart.}  

\subsection{Comparing the integrated optimization approach to the manual planning implemented in practice}
\label{sec:comparing}

In practice, SSTRPVST is solved sequentially by deciding the routing decisions of the sprayers, the refilling locations, and then the tanker's routing independently. More precisely, our partner company informed us that firstly they solve a TSP instance to find the shortest route connecting all nodes in the farm. The generated route is then partitioned into several routes based on the number of sprayers available to generate routing decisions for the sprayers. More formally, assume that the shortest route for a set of $\mathcal{N}^f$ nodes is written as an ordered sequence:
$$O: <0,i_1,i_2,i_3,i_4,i_5,...,i_n,0>$$

where $i_1$ is the first node visited, $i_2$ is the second node visited,..., and $i_n$ is the last node visited before returning to the depot. Route $O$ is then partitioned as follows assuming that two sprayers are available for spraying:

$$R_1: <0,i_1,i_3,i_5,...,i_n,0>$$
$$R_2: <0,i_2,i_4,i_6,...,i_{n-1},0>$$

where $R_1$ and $R_2$ are the ordered sequence of nodes visited by sprayer 1 and sprayer 2, respectively. The spraying quantity at each node is then assumed to be 10\% above the minimum required spraying quantity (i.e., $q_i=1.1*qMin_i$). At any node where the sprayer's tank level is less than the $q_i$ of the next node, a refilling is needed. With the refilling points established, the routing of the tanker can be established. If a waiting time at any node is expected, the service time is modified. The intuition behind this policy implemented in practice is that the sprayers will be working close to each other and hence waiting time is almost guaranteed to be zero. Furthermore, having the tanker follow a route generated by solving a TSP instance ensures that the tanker will not be driving around the farm inconsistently.  \\

In practice, the effectiveness of a spraying operation is measured by the total productivity of sprayers measured by the total service time of sprayers (i.e., the total time sprayers spend spraying) provided that the distances traveled by the tanker and the sprayers are minimized. The refilling times are undesired since they reduce the time available for the sprayer to spray. Clearly, the total service time of all sprayers and the routing time of sprayers and the tanker are affected by the policy used to generate the routing decisions of sprayers and therefore the routing of the tanker. Providing a policy that takes into account the integrated decisions of sprayers' routing, the tanker routing, and the service time of sprayers is expected to provide a more effective solution to a spraying operation as opposed to solving the same problem sequentially. To show it, we provide a computational experiment in this direction. First, we generate instances using real-world data for real-world farm area, namely a farm with 50 locations to be sprayed with a area of 2,000 acres. Then, we solve each instance using the matheuristic developed in our work and we compare the solution obtained using our matheuristic to the solution implemented in practice. To further provide some insights, we relax the requirement that the waiting time of the sprayers should be zero and we provide a model in which the waiting time of sprayers is allowed to gain some insights on what implications will be made if the waiting time was allowed. Therefore, we provide a comparison between three models, namely, the model presented in this paper where waiting time is not allowed, the policy implemented in practice, and finally a model where the waiting time is allowed. We generate 80 instances with different locations of nodes  to be sprayed and values of $qMin_i$ in a real-world large farm. Forty instances were solved assuming the availability of 3 sprayers and the other forty instances were solved assuming that 4 sprayers are available. \\

Table \ref{tab:modelsComparison} presents these average results for all models with 50 nodes. For each key performance measure, we report the value found using the optimization approach presented in this paper assuming that waiting time is allowed for sprayers, the approach used in practice (i.e., manual sequential planning), and the optimization approach presented in this paper. The former model represents a relaxation to model (\ref{eq:Objt})-(\ref{eq:FIntg}) and some trade-off is expected. However, it will be interesting to quantify this trade-off. \\

\begin{figure}[!h]
\centering
\includegraphics[width=14cm]{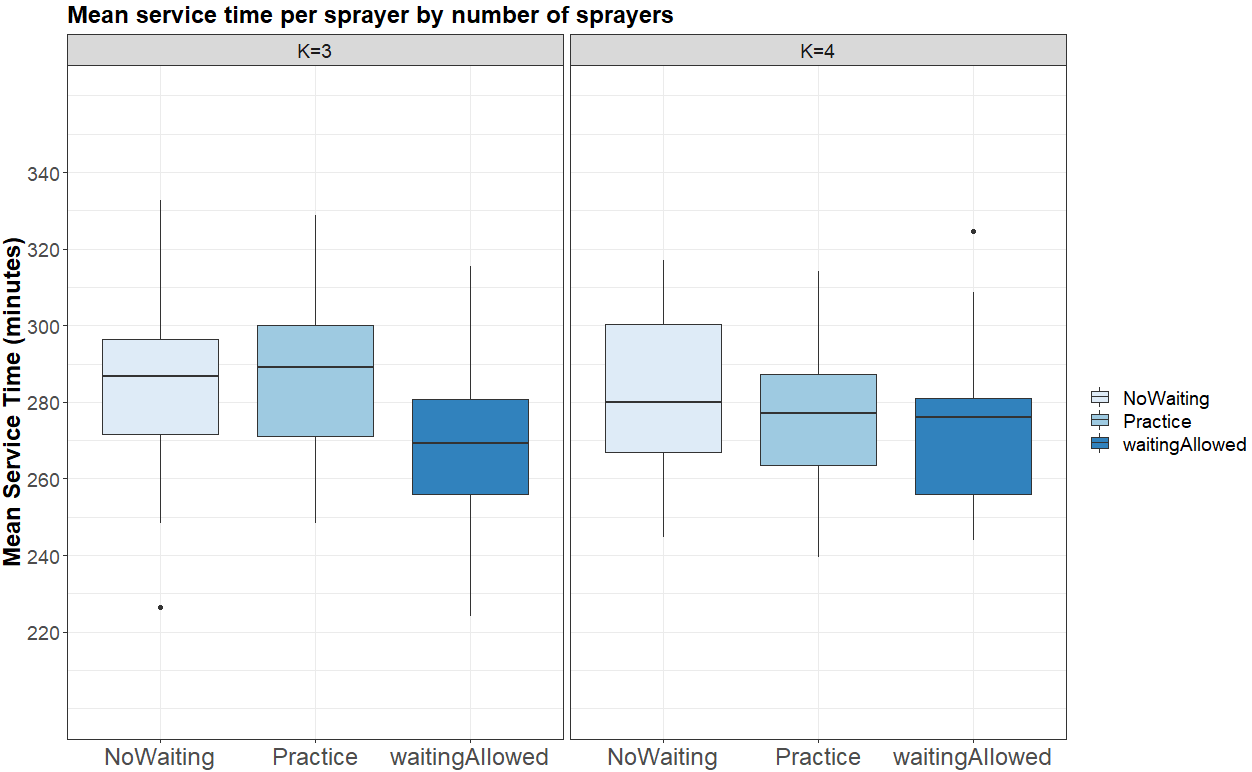}
\caption{Mean service time per sprayer under each model in minutes}
\label{fig:serModel}
\end{figure}

We first discuss the performance of each model on the \textit{the average service time per sprayer}. Having sprayers performing spraying operations is the essence of the whole spraying operation and it should not be surprising that farmers care about this performance measure. From Table \ref{tab:modelsComparison} we observe that models where waiting time is not allowed and the manual planning implemented in practice achieve similar values when it comes to service time. We note that the model that allows waiting time performs slightly worse than these two models. Specifically, the model implemented in practice achieves the highest service time per sprayer, followed by the optimization approach (the optimization approach yields a service time per sprayer that is $0.89\%$ less than the manual planning), then the model that allows waiting time achieves the lowest service level with a reduction of $3.42\%$ in service time compared to the manual planning. \textcolor{red}{Figure \ref{fig:serModel} shows the mean service time under each model by the number of sprayers in minutes. }\\

\begin{figure}[!h]
\centering
\includegraphics[width=14cm]{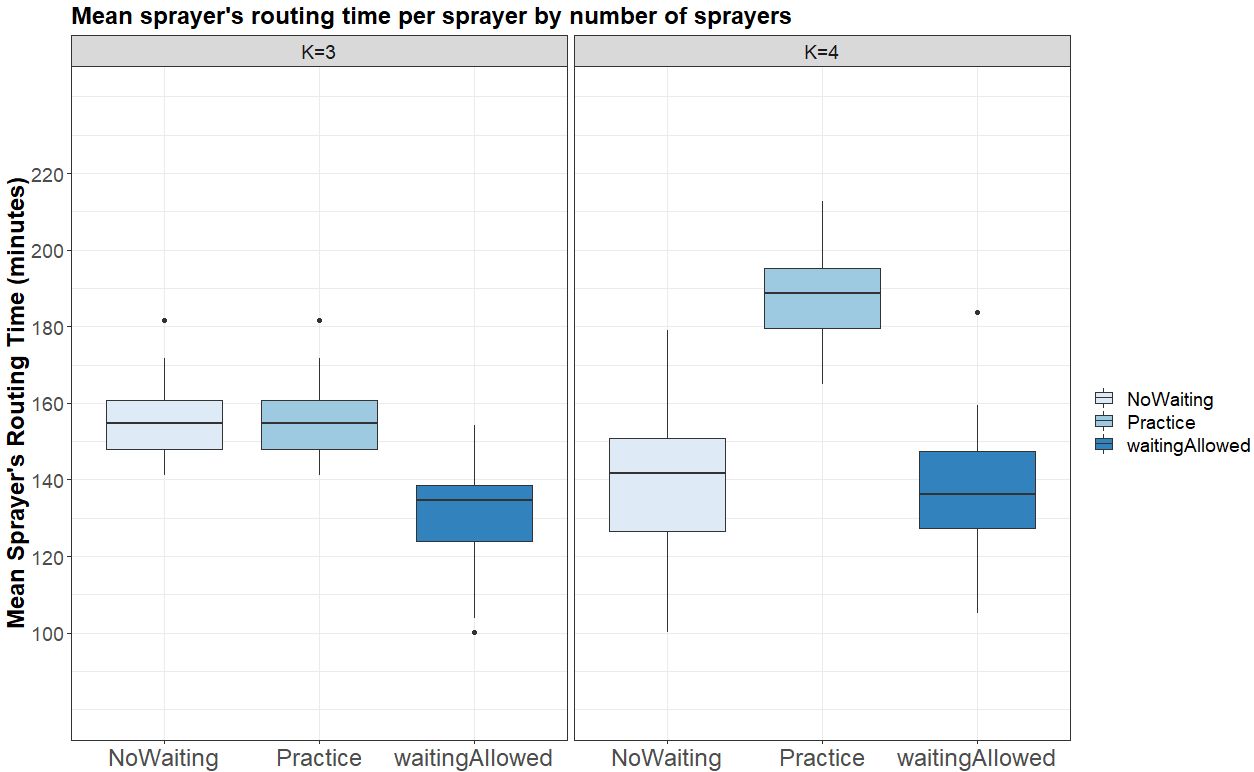}
\caption{Mean routing time per sprayer under each model in minutes}
\label{fig:routingModel}
\end{figure}

\begin{figure}[!h]
\centering
\includegraphics[width=14cm]{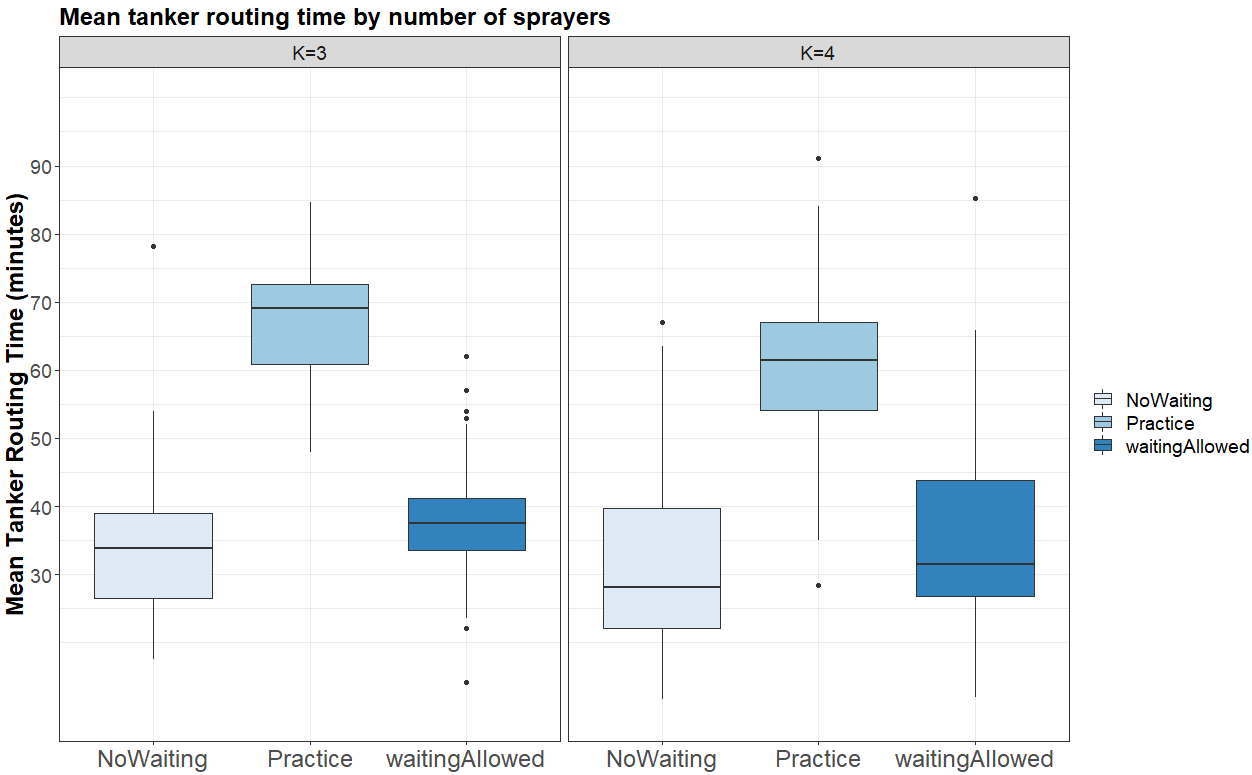}
\caption{Mean routing time of the tanker under each model in minutes}
\label{fig:tankerModel}
\end{figure}

Second, it is very critical for farmers to preserve healthy soil by minimizing the \textit{traveled distance of the sprayers and the tanker on the farm}. The model that achieves the best results in that performance measure is the model in which waiting time is allowed, followed by the optimization approach, and then the manual planning yields the highest traveling time for the sprayers and the tanker. These differences are rather significant, compared to the manual planning, the model allowing waiting time reduces the traveled distance per sprayer by $22.18\%$ and the optimization approach saves the traveled time per sprayer by $14.20\%$. These savings are rather significant and farmers will be satisfied if spraying can be performed while minimizing the traveled distance by such reductions. One clear explanation for this significant increase in the sprayer's routing time is that the manual planning does not provide a solution that integrates service time and routing decisions simultaneously. The same observation can be drawn for the tanker routing time. \textcolor{red}{Figures \ref{fig:routingModel} and \ref{fig:tankerModel} show the mean sprayer and tanker routing times under each model by the number of sprayers in minutes, respectively.}\\

Third, consider \textit{the number of re-fillings} needed. As shown in the second term in the objective function (\ref{eq:Objt}), we can see that minimizing the number of refills needed is an important part of the objective function since a sprayer is idle during the refilling. The model that achieves the best results in that performance measure is the model in which waiting time is allowed, followed by the optimization approach, and then the manual planning yields the highest traveling time for the tanker. It is very interesting that by reducing the service time by $3.41\%$ the model with waiting time can achieve a reduction of almost $40.0\%$ in the number of refills. \\

\begin{sidewaystable}[ph!]
\caption{Mean results of the service and routing times in minutes over the three models with 50 nodes instances}
\label{tab:modelsComparison}
\begin{tabular}{lllllllllllll}
 & \multicolumn{3}{c}{Service time} & \multicolumn{3}{c}{Sprayer routing time} & \multicolumn{3}{c}{Tanker routing time} & \multicolumn{3}{c}{Number of refills} \\
\multicolumn{1}{c}{$\vert \mathcal{K}\vert$} & \multicolumn{1}{c}{waiting all.} & \multicolumn{1}{c}{practice} & \multicolumn{1}{c}{no waiting} & \multicolumn{1}{c}{waiting all.} & \multicolumn{1}{c}{practice} & \multicolumn{1}{c}{no waiting} & \multicolumn{1}{c}{waiting all.} & \multicolumn{1}{c}{practice} & \multicolumn{1}{c}{no waiting} & \multicolumn{1}{c}{waiting all.} & \multicolumn{1}{c}{practice} & \multicolumn{1}{c}{no waiting} \\ \hline
3 & 270.56 & 286.61 & 281.92    & 130.73 & 155.21 & 155.66 & 38.02 & 66.96 & 34.42 & 5.78 & 7.71 & 5.00 \\
4 & 272.38 & 275.59 & 275.29    & 136.74 & 187.13 & 139.12 & 35.17 & 59.57 & 32.61 & 4.95 & 6.70 & 3.93 \\ \hline
Avg. & 271.47 & 281.10 & 278.61 & 133.73 & 171.85 & 147.45 & 37.01 & 63.26 & 33.52 & 5.37 & 7.20 & 4.46 \\ \hline
\end{tabular}
\end{sidewaystable}

\section{Conclusion}
\label{sec:conclude}

The SSTRPVST is a practical and highly complex problem. Considering multiple synchronization constraints and the variable service time in routing problems, to the best of our knowledge we were the first to address such complex problem. The SSTRPVST occurs in utilizing high-end spraying equipment in precision agriculture for which we have developed our solution approach. We built on an ALNS framework and enhance it with several innovative ideas, i.e. (a) new destroy operators, (b) an intensive local search, and (c) new neighborhoods exploiting the underlying problem structure. \\

We evaluated the matheuristic on a set of generic instances designed to represent different
farm layouts and different scenarios. Our computation results show that the developed enhancements add value to better solving the SSTRPVST. We generated insights into how different models present trade-offs in key performance measures. These insights can be used by companies implementing precision agriculture to decide which cost function to be used depending on their system and/or customer preferences. \\


\bibliographystyle{pomsref} 

 \let\oldbibliography\thebibliography
 \renewcommand{\thebibliography}[1]{%
 	\oldbibliography{#1}%
 	\baselineskip14pt 
 	\setlength{\itemsep}{10pt}
 }
\bibliography{ref1.bib} 

















\end{document}